\documentclass[12pt]{amsart}
\usepackage{amssymb}
\allowdisplaybreaks[1]

\theoremstyle{definition}

\numberwithin{equation}{section}

\newcommand{\T}{\mathbb T}
\newcommand{\R}{\mathbb R}

\renewcommand{\a}{\alpha}

\newcommand{\g}{\gamma}

\renewcommand{\v}{\varphi}
\renewcommand{\b}{\beta}
\newcommand{\s}{\sigma}
\renewcommand{\o}{\omega}

\renewcommand{\d}{\delta}
\newcommand{\D}{\Delta}

\renewcommand{\l}{\lambda}
\renewcommand{\L}{\Lambda}
\renewcommand{\O}{\Omega}
\renewcommand{\P}{\Phi}

\newcommand{\w}{\wedge}
\newcommand{\p}{\partial}
\newcommand{\dis}{\displaystyle}

\renewcommand{\P}{\mathcal P}

\newcommand\somme{\mathop {\oplus}\limits}

\newcommand{\union}{\mathop{\cup}}

\begin{document}

\title[An approach to the tangential Poisson cohomology]
{An approach to the tangential Poisson cohomology based on examples
in duals of Lie algebras}
\vskip3cm
\author[\small Angela Gammella]{Angela Gammella}
\address{
Universit\'e Libre de Bruxelles\\
Campus de la Plaine\\
CP 218, Boulevard du Triomphe\\
1050 Bruxelles, Belgique.}
\email{agammell@ulb.ac.be}


\vskip3cm

\maketitle

\

\

\

\begin{abstract}
We study the tangential Poisson cohomology
(TP-cohomology) of regular Poisson manifolds, first defined by Lichnerowicz using contravariant tensor
fields. We show that for a regular Poisson manifold
$M$, the TP-cohomology coincides 
with the leafwise de Rham (or $\check{\hbox{C}}$ech) cohomology of the symplectic foliation of $M$.
Its computation in various degrees leads to open, non-trivial problems. To get a better understanding of
these difficulties, we study explicitly many examples coming from nilpotent and 3-dimensional (real) Lie
algebras. For the latter, we compare the TP-cohomology and
the usual Poisson cohomology (P-cohomology).
\end{abstract}

\

\

\section{Introduction and motivation}\label{sec1}

This work fits into the study of deformation quantization for the dual $\frak{g}^\ast$ of a Lie algebra $\frak{g}$,
more exactly of star products on  $\frak{g}^\ast$ (or on some natural open subset $U$ of 
$\frak{g}^\ast$) which restrict nicely to
the coadjoint orbits contained in $\frak{g}^\ast$ (or $U$). 
Such star products are called tangential and for a given Lie algebra $\frak{g}$, they can notably be used 
to describe the harmonic analysis of the corresponding Lie group.\par
In general, tangential star products do not exist on the whole dual $\frak{g}^\ast$ (see
 [ACG] or [CGR]),
nevertheless we know there always exists such a star product
 on the dense subset $\O$ of maximal dimensional coadjoint orbits in $\frak{g}^\ast$ [Mas].
When studying all the possible classes of tangential star products on this set $\O$, 
we became interested in regular Poisson structures and especially in the TP-cohomology of regular Poisson manifolds.
(Indeed tangential star products are governed by the TP-cohomology; for instance 
their classification is described by the second TP-cohomology space.)
In this respect, it is worth mentioning Remarks 5 and 8 of the paper, which give an example
of how our work applies to the theory of tangential star products.\par
With this motivation from deformation theory, we present here
the result of our attempts 
to understand and clarify the TP-cohomology. We
organize the paper as follows.\par 
In Section 2 below, we prove that, for a regular Poisson manifold, the TP-cohomology 
is isomorphic to the leafwise de Rham (or $\check{\hbox{C}}$ech)
cohomology of the symplectic foliation. Thus it depends only on that foliation 
and not on the symplectic structure along the leaves.
This is a generalization of the fact that the P-cohomology of a symplectic manifold $M$ is just the de Rham cohomology
of $M$.    
We recall also some classical results of foliation theory computing the TP-cohomology of 
some particular regular Poisson
manifolds.  We compare these results with a theorem from [Va2] describing the P-cohomology for
some specific cases.\par
The remaining of the paper is devoted to examples arising from Lie algebras.
Indeed, each Lie algebra $\frak{g}$ gives rise to a natural regular Poisson manifold,
namely the union $\O$ of all maximal dimensional 
coadjoint orbits in the dual space $\frak{g}^\ast$.\par
In Section 3, we consider nilpotent Lie algebras $\frak{g}$.
For such Lie algebras, given a Jordan-H\"older basis $B$, $\O$ has a natural layering whose
first layer, say $V_B$, is known as the generic (dense) open subset of $\frak{g}^\ast$ associated to
$B$ ([ACG, Ver]). It is easy to see that the TP-cohomology of $V_{B}$ is trivial
in degree superior to zero. We prove here that the same is true for the
union $\dis\union_{B}V_{B}$, which is more canonical  
than $V_B$ since it does not depend on the choice of the basis $B$. However, 
$\dis\union_{B}V_{B}$ is sometimes strictly smaller than $\O$; this happens for instance in the case
of the filiform Lie algebras. We will also see, by studying in details the case of the filiform Lie algebra $\frak{g}_{4,1}$,
that the TP-cohomology of $\O$ can be essentially larger than the TP-cohomology of $\dis\union_{B} V_{B}$.
Such an example shows that the TP-cohomology of regular Poisson manifolds (and more generally the
leafwise de Rham cohomology of foliations) can be huge even if the leaves 
are cohomologically trivial.\par
Later, in Section 4, we consider an arbitrary 3-dimensional regular Poisson 
manifold $M$ and  we perform the inductive computations of [Va2] to describe the 
P-cohomology spaces of $M$. This enables us to observe the influence of the TP-cohomology on the P-cohomology:
the TP-cohomology spaces appear naturally in the decomposition of the 
P-cohomology spaces (see Proposition 6).
Then, we examine the TP-cohomology and the P-cohomology of the regular Poisson manifold $\O$ arising from 
any 3-dimensional Lie algebra. Some of these Lie algebras can be directly treated
with the help of Section 2, the others will require more attention. We conclude with some general remarks.

\

\section{Regular Poisson manifolds and foliation theory}\label{sec2}

\subsection{Basic definitions}

\

A Poisson manifold is a $C^{\infty}$ manifold $M$ equipped with a Poisson bracket $\{\,,\,\}$ 
{\sl i.e.} a bilinear skew-symmetric operation on 
$C^{\infty}(M)$ with values in $C^{\infty}(M)$, satisfying 
the Leibniz rule:\par
$$\{f,gh\}=\{f,g\}h+g\{f,h\}\quad\forall f,g,h\in C^{\infty}(M)$$
and the Jacobi identity:\par
$$\{\{f,g\},h\}+\{\{g,h\},f\}+\{\{h,f\},g\}=0\quad\forall f,g,h\in C^{\infty}(M).$$\par
For any manifold $M$, we denote by ${\mathcal V}^{\ast}(M)$ the graded space of skew-symmetric contravariant tensor fields 
and by $\O^{\ast}(M)$ the graded space of forms on $M$.\par
Let $M$ be a Poisson manifold.
Since the Poisson bracket is skew-symmetric and satisfies the Leibniz rule, there exists a unique tensor 
field $\L$ in ${\mathcal V}^2(M)$ such that 
$$\{f,g\}=\L(df,dg)\quad\forall f,g\in C^{\infty}(M).$$
This tensor field is usually called the Poisson bivector of $M$.\par
To express the 
Jacobi identity in terms of $\L$, we recall that the commutator bracket of vector
fields extends to the Schouten bracket, uniquely defined on ${\mathcal V}^{\ast}(M)$ by the relations:
\par\smallskip
(i)\, $[P,Q]=-(-1)^{(p-1)(q-1)}[Q,P]\quad \forall P\in {\mathcal V}^p(M),\forall Q\in {\mathcal V}^q(M)$\par\smallskip
(ii)\, For $P$ in ${\mathcal V}^p(M)$, $[P,.]$ is a derivation of degree $p-1$.
\par\smallskip\noindent
The Schouten bracket satisfies the graded Jacobi identity:
$$ [P, [Q,R]]= [[P,Q],R] +(-1)^{(p-1)(q-1)}[Q,[P,R]],$$
for  $P$ in ${\mathcal V}^p(M)$, $Q$ in ${\mathcal V}^q(M)$, $R$ in ${\mathcal V}^{\ast}(M)$, and thus
defines a graded Lie algebra structure on ${\mathcal V}^{\ast}(M)$ with the shifted grading:
$deg(S)=s-1$ if $S$ belongs to ${\mathcal V}^s(M)$. 
One can then check [Li1, Wei]
that the bracket on $C^{\infty}(M)$ given by $\L$ satisfies the Jacobi identity if and only if 
$[\L,\L]=0$ holds.\par 
In the sequel, we shall denote by $(M,\L)$ our Poisson manifold. If $f$ is a $C^{\infty}$ function on $M$, 
we call
Hamitonian vector field of $f$ the vector field corresponding to the derivation $\{f,.\}$.
With $(M,\L)$ is associated a bundle map:
\begin{align*}
\#:T^{\ast}M&\longrightarrow TM\\
\a&\longmapsto \a^{\#}
\end{align*}
defined by $$\a^{\#}(\b)=\L(\a,\b)$$ for any $\a,\b$ in $T^{\ast}_xM$.
Finally, the rank of $M$ at a point $x$ is by definition the rank of the linear mapping 
$\#_x:T_x^{\ast}M\rightarrow T_xM$. 
If it is constant, $M$ is said to be regular. In particular, if it is 
everywhere equal to the dimension of $M$, $\#$ is
 an isomorphism and $M$ is a symplectic manifold whose
symplectic structure $\o$ is given by $\o=\#^{-1} (\L)$.\smallbreak
We are now ready to define the Poisson cohomology of the Poisson manifold $(M,\L)$. 

\newtheorem{Def}{Definition}
\begin{Def}
Let $\s:{\mathcal V}^{\ast}(M)\rightarrow {\mathcal V}^{\ast+1}(M)$ be the operator given by
$$\s=[\L,.].$$
Due to the graded Jacobi identity for $[\,,\,]$, 
$\s$ is a coboundary operator ({\sl i.e.} $\s^2=0$).
The complex $({\mathcal V}^{\ast}(M),\s)$ is called the Poisson complex of $M$  and the corresponding
cohomology 
$H^{\ast}_{\L}(M)$ is the P-cohomology of $M$. 
\end{Def}

See [Hue] for an algebraic definition of the P-cohomology and [APP, CW] 
for some more recent results about the P-cohomology.\par
The interpretation of the first few P-cohomology spaces is well-known. Indeed,
$H^0_{\L}(M)$ is the space $I(M)$ of Casimir functions over $M$ {\sl i.e.} those whose Hamiltonian vector 
fields are trivial;
$H^1_{\L}(M)$ consists of infinitesimal Poisson automorphisms of $M$ (Poisson vector fields)
 modulo inner automorphisms (Hamiltonian vector fields);
$H^2_{\L}(M)$ classifies (modulo the trivial deformations)
the formal deformations of the Poisson structure $\L$ (with the form $\L + t \a_1+ t^2\a_2+\ldots$)
and finally $H^3_{\L}(M)$ houses the obstructions to extend a formal deformation  
from one step (in powers of $t$) to the next.\par
Let us just recall that the equivalence classes of star products on $(M, \L)$ are in one to one
correspondence with the equivalence classes of formal deformations of $\L$
(see [Kon] for more details).\par
Note also that the Poisson bracket gives rise to a bracket $\{\,,\,\}$ on $\O^1(M)$, which is
the unique extension of the bracket given by
$\{df,dg\}=d\{f,g\}$
such that
$$\{\a,f \b\}=f\{\a,\b\}+(\a^{\#} f) \b \quad \forall f\in C^{\infty}(M)\quad\forall \a,\b\in \O^1(M).$$
This bracket is defined by
$$\{\a,\b\}=L_{\a^{\#}}(\b)-L_{\b^{\#}}(\a)-d(\L(\a,\b))\quad\forall \a,\b\in \O^1(M)$$
and one can prove (see [Va2] p.44) 
\begin{align*}
&\s Q({\a}_0,...,{\a}k)=\sum_{i=0}^k (-1)^i  \a_i^{\#}(Q(\a_0,\ldots,\hat{\a_i},\ldots,\a_k))+\cr
&\sum_{i<j} (-1)^{i+j} Q(\{\a_i,\a_j\},\a_0,\ldots,\hat{\a_i},\ldots,\hat{\a_j},\ldots,\a_k)
\end{align*}
where $Q$ is in ${\mathcal V}^k(M)$ and the $\a_i$ are 1-forms on $M$.\par\noindent
The latter expression can be used ([Va2]) to see that the natural extension $\tilde{\#}$ of $\#$ to forms:
$$\tilde{\#}\l(\a_0,\ldots,\a_{q-1})=(-1)^q \l(\a_0^{\#},\ldots,\a_{q-1}^{\#})$$
intertwines $\s$ and the de Rham differential $d$ and
thus induces a natural  homomorphism from $H^{\ast}_{DR}(M)$ to $H_{\L}^{\ast}(M)$.
This homomorphism is trivially an isomorphism in the symplectic case
(see also [Kos] or [Li1]).\par

\medbreak
Some preparatory material related to a foliation is now needed. Let  $(M,{\mathcal F})$
be an arbitrary foliated manifold and denote by $T{\mathcal F}$ the tangent
bundle of ${\mathcal F}$.
As in [DH] or [Li2], one can choose a transversal distribution $\nu{\mathcal F}$ such that 
$$TM=T{\mathcal F}\oplus \nu{\mathcal F}\quad\hbox{and}\quad T^{\ast}M=T^{\ast}{\mathcal F}
\oplus\nu^{\ast}{\mathcal F}.$$
These decompositions 
induce a bigrading of the space ${\mathcal V}^{\ast}(M)$ of contravariant
tensor fields  and of the space 
$\O^{\ast}(M)$ of forms on $M$, namely
$${\mathcal V}^{\ast}(M)=\somme_{p,q}{\mathcal V}_{p,q}(M)\quad\hbox{and}\quad
\O^{\ast}(M)=\somme_{p,q}\O_{p,q}(M)$$
where
$\dis{\mathcal V}_{p,q}(M)$ (resp. $\O_{p,q}(M)$)
denotes the space of sections of the bundle    $\w^q(T{\mathcal F})\otimes \w^p(\nu{\mathcal F})$
(resp. $\w^q(T^{\ast}{\mathcal F})\otimes \w^p({\nu}^{\ast}{\mathcal F})).$
\par\noindent
Elements of ${\mathcal V}_{p,q}(M)$ and $\O_{p,q}(M)$ are said to be of type $(p,q)$.
Moreover, an operator will be homogeneous of type $(a,b)$ if it sends elements of type $(p,q)$ to elements of type
$(p+a,q+b)$. 
We recall that the de Rham differential $d$ can be decomposed into $d=d'+d''+d_{2,-1}$ where
$d'$ is of type $(1,0)$, $d''$ denotes the leafwise de Rham differential of the foliation ${\mathcal F}$
and is of type $(0,1)$, and $d_{2,-1}$ is of type $(2,-1)$.\medbreak
Assume now that $(M,\L)$ is a regular Poisson manifold and denote by ${\mathcal F}$ the symplectic foliation 
of $M$. As above, one can choose a transversal distribution $\nu{\mathcal F}$ for $M$.
It was shown in [Va2] that with respect to a given choice of $\nu{\mathcal F}$, the coboundary operator
$\s$, introduced in Definition 1,
has a well defined decomposition $\s=\s'+\s''$ where $\s'$ is of type $(-1,2)$ and $\s''$ is of type $(0,1)$.\par\noindent
On the other hand, Lichnerowicz has shown in [Li2] that one gets a consistent theory by restricting 
the P-cohomology complex
$({\mathcal V}^{\ast}(M),\s)$ to tangential multivector fields. The resulting cohomology is known as the TP-cohomology
of the regular Poisson manifold $(M,\L)$. 
In fact, the same cohomology can be defined by using the transversal distribution $\nu{\mathcal F}$ and
the types of the tensor fields. Indeed, we have
\begin{Def} The TP-cohomology complex of the regular Poisson manifold $(M,\L)$ is 
$\dis\somme_{q}{\mathcal V}_{0,q}(M)$ with the coboundary operator $\s''$ 
 and 
$$H_{\L,tan}^{q}(M)=\frac{ Ker (\s'':{\mathcal V}_{0,q}(M)\rightarrow {\mathcal V}_{0,q+1}(M))}
{Im(\s'':{\mathcal V}_{0,q-1}(M)\rightarrow {\mathcal V}_{0,q}(M))}$$
is the {\sl qth} TP-cohomology space of $(M,\L)$.\end{Def}\par
It is clear that
$H^0_{\L,tan}(M)=H^0_{\L}(M)$.
Moreover, as mentioned in the introduction, the TP-cohomology plays an important role in the theory of 
tangential star products. Indeed,
the derivations of a given tangential star product on $M$, modulo inner derivations, are
parameterized by sequences of elements in $H^1_{\L,tan}(M)$; similarly  equivalences of tangential star products on $M$
are classified at each step by $H^2_{\L,tan}(M)$ and  finally the obstructions to
construct such a star product are localized in $H^3_{\L,tan}(M)$
(this last point could be omitted since a tangential deformation always exists on $M$ [Mas]).
\vskip0,25cm
\subsection{Leafwise de Rham cohomology}

\

In this paragraph, we want to prove that the TP-cohomology of a regular Poisson manifold $(M,\L)$ 
is isomorphic to the leafwise de Rham cohomology of the symplectic foliation 
and therefore does not depend on the symplectic 
structure along the leaves.\par
To this end, we consider first the general case of a foliated manifold $(M,{\mathcal F})$.
For all $p$, we denote by $\P^p({\mathcal F})$ the sheaf of (germs of) projectable $p$-forms ({\sl i.e.}
those induced by forms on the space of leaves). In particular, 
$\P^0({\mathcal F})$ is the sheaf of the germs of functions on $M$ that are constant along the 
leaves of ${\mathcal F}$ (see [Va2]).
Let us fix a transversal distribution $\nu{\mathcal F}$ and consider the sheaf cohomology of 
$\P^p({\mathcal F})$, that is
$$H^q(M,\P^p({\mathcal F}))=\frac{ Ker (d'':\O_{p,q}(M)\rightarrow \O_{p,q+1}(M))}
{Im (d'':\O_{p,q-1}(M)\rightarrow \O_{p,q}(M))}.$$
A change of $\nu{\mathcal F}$ leads to an isomorphism in the corresponding cohomology spaces. To see this, let
us introduce more notations. 
Let $N{\mathcal F}$ be the normal bundle of ${\mathcal F}$, that is $N{\mathcal F}=TM/ T{\mathcal F}$. Denote by 
$\O_{[p]}^q({\mathcal F})$ the space of sections of the bundle
$\w^q(T^{\ast}{\mathcal F})\otimes \w^p(N^{\ast}{\mathcal F})$.
The elements of $\O_{[p]}^0({\mathcal F})$ are usually called normal forms and
those of $\O_{[p]}^q({\mathcal F})$ are the tangential $q$-forms with values
in the normal $p$-forms. The Lie algebra of tangential vector fields acts
naturally (by Lie derivative) on the normal forms and the leafwise de Rham differential
$d_{\mathcal F}$ acts on $\O_{[p]}^q({\mathcal F})$ in the usual way:
\begin{align*}
d_{\mathcal
F}(\o_{x})(X_0,\ldots,X_q)=&\dis\sum_{i=0}^q
(-1)^iL_{X_i}(\o_{x}(X_0,\ldots,\hat{X_i},\ldots, X_q))+\cr
&\dis\sum_{i<j}(-1)^{i+j}\o_{x}([X_i,X_j],X_0,\ldots, \hat{X_i},\ldots, \hat{X_j},\ldots, X_q),
\end{align*}
if $x$ is in $M$, $\o$ in $\O^q_{[p]}({\mathcal F})$ and the $X_i$ in $T_{x}{\mathcal F}$. 
Since $d_{\mathcal F}^2=0$, we obtain a complex
 $(\dis\somme_{q}\O^{q}_{[p]}({\mathcal F}),d_{\mathcal F})$ whose cohomology is denoted hereafter by 
$H_{[p]}^{\ast}({\mathcal F})$.
\par\noindent
With these notations, we can prove 
\newtheorem{Lem}{Lemma}
\begin{Lem} For all $p$, 
the complexes
$(\dis\somme_{q}\O_{p,q}(M),d'')$ and $(\dis\somme_{q}\O^{q}_{[p]}({\mathcal F}),d_{\mathcal F})$ are isomorphic. In particular, 
$H^{\ast}(M,\P^p({\mathcal F}))$ and $H_{[p]}^{\ast}({\mathcal F})$ coincide and are independent of the choice of 
$\nu{\mathcal F}$.
\end{Lem}\par
{\sl Proof.}
For each $x$ in $M$, let $\psi_{x}:N_{x}{\mathcal F}\rightarrow \nu_{x}{\mathcal F}$ be the natural isomorphism
of vector spaces. That is,
if $Y$ is in $T_{x}M$, then $\psi_{x}(\tilde{Y})=\pi(Y)$ where
$\tilde{Y}$ denotes the class of $Y$ in $N_{x}{\mathcal F}$ and $\pi$ is the projection
$\pi:TM\rightarrow \nu{\mathcal F}$.\par\noindent
The mapping $f:\O_{p,q}(M)\rightarrow \O_{[p]}^q(M)$ given by
$$f(\o_{x})(X_1,\ldots,X_q)(\tilde{Y_1},\ldots,\tilde{Y_p})=
\o_{x}(X_1,\ldots,X_q,\psi_{x}(\tilde{Y_1}),\ldots,\psi_{x}(\tilde{Y_p})),$$
here $\o$ is in $\O_{p,q}(M)$, the $X_i$ in $T_{x}{\mathcal F}$, and the  $Y_j$ in $T_{x}M$,
is clearly bijective.
In addition, we have
\begin{align*}
&d_{\mathcal F}f(\o_{x})(X_0,\ldots,X_q)(\tilde{Y_1},\ldots,\tilde{Y_p})\cr
&=\sum_{i=0}^q(-1)^i
X_i\Bigl(f(\o_{x})(X_0,\ldots,\hat{X_i},\ldots,X_q)(\tilde{Y_1},\ldots,\tilde{Y_p})\Bigr)\cr
&-\sum_{i=0}^q\sum_{l=1}^p(-1)^{i+l-1}f(\o_{x})(X_0,\ldots,\hat{X_i},\ldots,X_q)
(\widetilde{[X_i,Y_l]},\tilde{Y_1},\ldots,\hat{Y_l},\ldots,\tilde{Y_p})\cr
 &+\sum_{0\leq i<j\leq
q}(-1)^{i+j} f(\o_{x})([X_i,X_j],X_0,\ldots,\hat{X_i},\ldots,\hat{X_j},\ldots X_q)
(\tilde{Y_1},\ldots,\tilde{Y_p})\cr
 &=f(d''\o_{x})(X_0,\ldots,X_q)(\tilde{Y_1},\ldots,\tilde{Y_p}).
\end{align*}
Thus, $d_{\mathcal F}\circ f=f\circ d''$. This ends the proof.\par\noindent
Remark that  $H^q(M,\P^0({\mathcal F}))$ (or $H_{[0]}^q({\mathcal F})$) is nothing else 
but the leafwise de Rham cohomology of the foliation ${\mathcal F}$.\medskip
\par
Let us now restrict ourselves to the case of a regular Poisson manifold $(M,\L)$.
Recall that $I(M)$ denotes the space of Casimir functions over $M$, in other words the space of 
smooth functions on $M$ that are constant along the leaves of the symplectic foliation.
Then we have

\newtheorem{Th}{Theorem}
\begin{Th}
Let $(M,\L)$ be a regular Poisson manifold, ${\mathcal F}$ the symplectic foliation of $M$ and
$\nu{\mathcal F}$ a transversal distribution for $M$.
Then, $(\dis\somme_{q} {\mathcal V}_{0,q}(M),\s'')$ and $(\dis\somme_{q} \O_{0,q}(M),d'')$ are isomorphic as complexes 
of $I(M)$-modules. In particular, for all $q$, $H^q_{\L,tan}(M)$ and $H^q(M,\P^0({\mathcal F}))$
are isomorphic $I(M)$-modules.
\end{Th}
\noindent\par
{\sl Proof.} 
It is not difficult to prove that the natural extension $\tilde{\#}$ 
of $\#$ realizes an $I(M)$-modules isomorphism between  
 $\O_{0,q}(M)$ and ${\mathcal V}_{0,q}(M)$. This isomorphism satisfies $$\s''\circ
\tilde{\#}=-\tilde{\#}\circ d''.$$  The result follows.
\par\medbreak
\subsection{$\check{\hbox{C}}$ech cohomology.}

\

Let $(M,{\mathcal F})$ be a foliated manifold. We shall say that a locally finite covering  
${\mathcal U}=(U_i)$ of $M$ is a good covering if 
for all $q>0$, all $k$ and all $i_1,...,i_{k}$,  
$$H^q(U_{i_1,...,i_{k}},\P^0({\mathcal F}_{|U_{i_1,...,i_{k}}}))=\{0\},$$
where $U_{i_1,...,i_k}=U_{i_1}\cap\ldots\cap U_{i_k}$ and 
${\mathcal F}_{|U_{i_1,...,i_{k}}}$ denotes the foliation induced by ${\mathcal F}$ in $U_{i_1,...,i_k}$
({\sl i.e.} if the leaves of ${\mathcal F}$ are  noted ${\mathcal L}_{\a}$, the leaves of 
${\mathcal F}_{|U_{i_1,...,i_{k}}}$ are the connected components of the intersections
${\mathcal L}_{\a}\cap U_{i_1,...,i_{k}}$).\par

We recall that, for each foliated manifold $(M,{\mathcal F})$, there exists affine connections on $M$,
which are torsion free and adapted to ${\mathcal F}$ in the sense of [Li2].
Let us now prove the existence of good coverings.

\begin{Lem}
Let $(M,{\mathcal F})$ be a foliated manifold and $\Gamma$ an affine connection on $M$, which is torsion
free and adapted to ${\mathcal F}$.
Then, every atlas $\{U_i,\v_i\}$ on $M$ of distinguished charts has 
an open refinement $\{V_l,{\psi}_l\}$ such that\par\noindent
(1) each $V_l$ has compact closure\par\noindent
(2) $(V_l)$ is locally finite and is a good covering of $M$.\par\noindent
 
\end{Lem}\par
{\sl Proof.}
By taking an open refinement if necessary, we may assume that $(U_i)$ is 
locally finite and that each $U_i$ has compact closure. Let $(U_i')$ be an open refinement of $(U_i)$
(with the same index set) such that $\overline{U_i'}\subset U_i$ for all $i$, and note 
 $\v_i'={\v_i}_{|_{U_i'}}$. 
For each $x$ in $M$, let $W_{x}$ be a normal neighborhood of $x$, which is small enough to
satisfy the following properties:\par\noindent
- for each point $a$ in $W_x$, there exists a normal neighborhood $N_a$ of $0$ in $T_aM$ such that
$\exp:N_a\rightarrow W_x$ is a diffeomorphism;\par\noindent
- $W_x$ is geodesically convex;\par\noindent
- $W_x$ is contained in some $U_i'$.\par\noindent
Note $\v_x={\v_i'}_{|_{W_x}}$ and, 
for each $k$, set $$B_k=\{(W_x,\v_x): W_x\cap \overline{U_k'}\not= \emptyset\}.$$
Since $\overline{U_k'}$ is compact, there exists a finite subfamily $B_k'$ of $B_k$, which covers
$\overline{U_k'}$.
Then the family ${\mathcal B}=\dis\union_{k} B_k'$ is an open refinement $\{V_l,\psi_l\}$ of $\{U_i,\v_i\}$, which satisfies
(1) and (2).
In fact, it is clear by construction that each $V_l$ has compact closure and that
$(V_l)$ is locally finite. To show that $(V_l)$ is also a good covering of $M$, we shall now prove
that each $V_l$ has geodesically convex plaques.
Let $y$ be in some $V_l$ and denote by 
$P_y$ the plaque of ${\mathcal F}_{|_{V_l}}$ containing $y$.
Take two points $a$ and $b$ in $P_y$.
By assumption,
there exists a normal neighborhood $N_a$ of $0$ in $T_a M$ such that
$\exp:N_a\rightarrow V_l$ is a diffeomorphism.\par
a) First, we shall prove that $\exp(N_a\cap T_a{\mathcal F})=P_a$.
Let $Y$ be in $N_a\cap T_a{\mathcal F}$. Denote by $\tau$ the geodesic of $V_l$ with the initial condition $(a, Y)$:
$$\tau(0)=a\quad\hbox{and}\quad\dot{\tau}(0)=\tau_{\ast_{0}}(\frac{d}{dt})=Y.$$
This curve $\tau$ is at least defined for $0\leq t\leq 1$.
We express it in the distinguished chart $(V_l,\psi_l)$ as follows
$$x(t):=\psi_l(\tau(t))=(x^i(t),x^u(t))$$
where $i,j\ldots=1,\ldots,r$ (resp.
$u,v,\ldots=1,\ldots,s$) denote the tangential (resp. transverse) indexes.
Since $\tau$ is a geodesic, it satisfies

\begin{align*}
&\dis\frac{d^2 x^i}{dt^2}=-\dis\sum_{1\leq J,K\leq r+s} \Gamma_{J K}^i(x) 
\frac{dx^J}{dt}\frac{dx^K}{dt} \quad\hbox{for}\quad i=1,\ldots,r\\
&\dis\frac{d^2 x^u}{dt^2}=-\dis\sum_{1\leq J,K\leq r+s} \Gamma_{J K}^u(x)
\frac{dx^J}{dt}\frac{dx^K}{dt} \quad\hbox{for}\quad u=1,\ldots,s.
\end{align*}

Moreover, since $\Gamma$ is torsion free and adapted to ${\mathcal F}$, we have
$$\Gamma_{iA}^u=\Gamma_{Ai}^u=0 \quad\forall 1\leq i\leq r,\,1\leq u\leq s, \,1\leq A\leq r+s.$$ 
It follows that
$$\frac{d^2 x^u}{dt^2}=-\dis\sum_{1\leq v,w\leq s} \Gamma_{v w}^u(x) \frac{dx^v}{dt}
\frac{dx^w}{dt}\quad
\hbox{for}\quad u=1,\ldots,s.$$
Let us introduce the notation
$$f(t)=(x^i(t))_{1\leq i\leq r},\quad g(t)=(x^u(t))_{1\leq u\leq s}.$$
Then, the above system can be reduced to two ordinary differential equations of the form
\begin{align*}
(1)\quad f''(t)&=F(f'(t),f(t),g'(t),g(t))\\
(2)\quad g''(t)&=G(g'(t),g(t),f(t)).
\end{align*}
Since $\dot{\tau}(0)=Y$ is in $T_a{\mathcal F}$, we shall have $g'(0)=0$.
Now, for fixed $f$ and with the initial conditions
$g(0)$ and $g'(0)=0$, (2) has a unique solution namely $g=cst=g(0)$.
Let us denote by
$f_0(t)=(a^i(t))_{1\leq i\leq r}$ the unique solution of (1) when $g=cst$ and with $f(0)$ and $f'(0)$
as initial conditions.
Then, we have
$(x^i(t),x^u(t))=(a^i(t),x^u(0))$ for all $0\leq t\leq 1.$
In particular, 
$\tau(1)=\exp(Y)$ belongs to $P_a$.
We get thus the inclusion $\exp(N_a\cap T_a{\mathcal F})\subset P_a.$
The equality comes from the fact that $\exp(N_a\cap T_a{\mathcal F})$ is both open and closed in $P_a$, and that $P_a$
is connected.
\par
b) Now, let $\g$ be the unique minimizing geodesic of $V_l$, joining $a$ and $b$. We may write
$$\g:[0,1]\rightarrow V_l,\quad \g(t)=exp(tX)$$
where $\g(0)=a$, $\g(1)=b$ and $X$ is in $N_a$.
Since $\g(1)=b=\exp(X)$ is in $P_a(=P_y)$ and using the equality proved in a), we see that
$X$ is in fact in $N_a\cap T_a{\mathcal F}$. Thus,
$\g$ lies entirely in $P_y$.\par

We have proved that $V_l$ has geodesically convex (hence contractible) 
plaques. It is of course the same for every finite intersection 
$V_{l_1,\ldots,l_k}$. Therefore, each foliation ${\mathcal F}_{|_{V_{l_1,\ldots,l_k}}}$
is a product foliation by contractible leaves. Following a result of [Va2],
we shall mention just in the next paragraph (Theorem 3), this means that
$H^q(V_{l_1,\ldots, l_k},\P^0({\mathcal F}_{|_{V_{l_1,\ldots,l_k}}}))=\{0\}$ for all $q>0$. Lemma 2 is proved.
\par\medbreak

Now, for any sheaf $A$ on $M$ (in particular for $\P^0({\mathcal F})$), 
we shall denote by
$\check{C}^{k}({\mathcal U},A)$ the space of $k$-$\check{\hbox{C}}$ech cochains of $A$ with respect to a 
covering ${\mathcal U}$ of $M$, by $\check{\d}$ the $\check{\hbox{C}}$ech coboundary and 
$\check{\hbox{H}}^{\ast}({\mathcal U},A)$ the cohomology
corresponding to the $\check{\hbox{C}}$ech complex ($\check{C}^{\ast}({\mathcal U},A),\check{\d})$.\par
The purpose of the following proposition is to prove that the leafwise de Rham cohomology of 
a foliation ${\mathcal F}$ 
coincides with the $\check{\hbox{C}}$ech cohomology of the sheaf $\P^0({\mathcal F})$. This can be convenient 
to calculate the TP-cohomology (see later in Section 3).
 
\newtheorem{Pro}{Proposition}
\begin{Pro}
Let $(M,{\mathcal F})$ be a foliated manifold and
 $\nu{\mathcal F}$ a transversal distribution for $M$.
Let also ${\mathcal U}=(U_i)$ be a good covering of $M$ and $(h_i)$  a partition of unity subordinate to ${\mathcal U}$.
For a $k$-$\check{\hbox{C}}$ech cocycle $c$, let $\o_c$ be the d''-closed $k$-form on $M$ defined by
$$\dis{\o_c}_{|U_i}=(-1)^{\frac{k(k+1)}{2}}\sum_{i_1,...,i_k} c_{{i_1,...i_k},i} d''h_{i_1}\w\ldots
d''h_{i_k}.$$ Then, the homomorphism $\v: \check{H}^k({\mathcal U},\P^0({\mathcal F}))\rightarrow
H^k(M,\P^0({\mathcal F}))$ mapping the cohomology class $[c]$ to the cohomology class $[\o_c]$ is an
isomorphism for all $k$.
\end{Pro}\par
{\sl Proof.} 
The result can be proved in the same way as Brylinski did in [Bry] to show that the $\check{\hbox{C}}$ech cohomology of 
the constant sheaf
${\R}_{M}$ on any manifold $M$ coincides with the de Rham cohomology of $M$.
The key point of the proof is to consider the $\check{\hbox{C}}$ech double complex 
$K^{\ast\ast}=\check{C}^{\ast}({\mathcal U},\underline{\O_{0,{\ast}}(M)})$ where underlining denote
sheaves of germs. The total cohomology of $K^{\ast\ast}$, also called $\check{\hbox{C}}$ech
hypercohomology, is by definition the cohomology corresponding to the complex
$(K^{\ast},D)$ where $K^n=\somme_{p+q=n} K^{p,q}$  and $D=\check{\d}+(-1)^p d''$ in degree $(p,q)$.
It is known ([Bry] p.28) that the natural spectral sequences associated to $K^{\ast}$ 
lead to a canonical isomorphism between $\check{H}^k({\mathcal U},\P^0({\mathcal F}))$ and  $H^k(M,\P^0({\mathcal F}))$.
To prove this isomorphism is just $\v$, one needs first to see
$c$ and $\o_c$ as elements of $\check{C}^k({\mathcal U},\underline{\O_{0,0}(M)})$ and
$\check{C}^0({\mathcal U},\underline{\O_{0,k}(M)})$ respectively. Then one can show, as in [Bry] p.45,
that
$[c]$ and $[\o_c]$ correspond to each other from the point of view of hypercohomology. The proposition
follows directly.

\vskip0,25cm
\subsection{Review of some classical results}

\

We begin this paragraph by mentioning two well-known results of foliation theory.
These results, which can be found in [DH] and [Va2] respectively, 
provide the computation of the TP-cohomology in some particular cases.

\begin{Th}
Let $(M,{\mathcal F})$ be a foliated manifold and $r$  some integer.
Assume that the foliation ${\mathcal F}$ is given by a submersion 
$\Pi:M\rightarrow B$ ($B$ being a Hausdorff manifold), and that 
any leaf $L$ of ${\mathcal F}$ is connected and satisfies
$H_{DR}^q(L)=\{0\}$ for all $0<q\leq r$. Then,
\[ H^q(M,\P^0({\mathcal F}))=
\left\{\begin{array}{ll}
C^{\infty}(B) & \hbox{if $\,q=0$}\\
\{0\}& \hbox{if $\,0<q\leq r$.}
\end{array}
\right.\]
\end{Th}

\begin{Th}
Let $L$ and $R$ be two smooth manifolds and  ${\mathcal F}$  the foliation of $M=L\times R$ by the leaves $L\times\{x\}$ 
where $x$ is in $R$. Assume 
that $L$ has finite Betti numbers. Then,
$$H^q(M,\P^p({\mathcal F}))=H_{DR}^q(L)\otimes \O^p(R).$$
\end{Th}

Next, in Theorem 4, we recall an important result from [Va1, Va2]. This result will be used in Section 4 when we shall
compute the P-cohomology spaces associated to 3-dimensional Lie algebras.

\begin{Th}
Let $M=S\times R$ be a regular Poisson manifold 
whose regular Poisson structure $\L$ is transversally constant with respect to the transversal distribution
$\nu{\mathcal F}=T R$ ({\sl i.e.} the symplectic foliation ${\mathcal F}$ of $M=S\times R$ is defined by a fixed symplectic
structure of $S$).
Suppose that $S$ has finite Betti numbers. Then,
$$H_{\L}^q(M)\simeq\somme_{0\leq k\leq q} H_{DR}^k(S)\otimes \O^{q-k}(R).$$
\end{Th}

\newtheorem{Rem}{Remark}
\begin{Rem}
One can use Theorem 4 to show that, with constrast to the TP-cohomology, the P-cohomology not only 
depends on the symplectic foliation but also on the symplectic structure along the leaves.
Indeed, let $M$ be ${\bf S}^2\times {\R}^{\ast}_{+}$ and denote by 
$\o$ the standard symplectic structure on the unit sphere ${\bf S}^2$.
If $M$ is endowed with the regular Poisson structure defined by the same symplectic structure $\o$ on each leaf, 
then the P-cohomology of $M$ is given by Theorem 4.
But, if the same manifold $M$ is viewed as $\frak{su}(2)^{\ast}\backslash\{0\}$ with its usual
linear Poisson structure, then
each leaf ${\bf S}^2\times \{ t \}$ ($ t\in {\R}^{\ast}_{+}$) has a different symplectic structure,
namely
$t\o$, and Theorem 4 is no more valid for $M$ (we will see the actual computation of 
$\frak{su}(2)^{\ast}\backslash\{0\}$ in ¤4.7).
\end{Rem}
\smallskip\par
Apart from the specific cases of Theorems 2 and 3, the task of computing the 
TP-cohomology still remains unsolved. To better understand the TP-cohomology of general 
regular Poisson manifolds and to make some comparison between the TP-cohomology and the P-cohomology, 
we devote the next two sections to a large number of explicit computations related to Lie algebras.
More precisely, the regular Poisson manifolds we shall consider in the rest of the paper
are Poisson submanifolds of the union $\O$ of all maximal dimensional coadjoint orbits 
in the dual of a given Lie algebra.
That makes sense since the dual $\frak{g}^\ast$ of any Lie algebra $\frak{g}$ can be endowed
with a natural Poisson structure, the well-known Lie Poisson structure [We]; the leaves of the symplectic foliation of
$\frak{g}^\ast$ being exactly the coadjoint orbits.
\vskip0,5cm

\

\section{The nilpotent case}\label{sec3}

Suppose that $\frak{g}$ is an $m$-dimensional nilpotent Lie algebra. Denote by $\frak{g}^{\ast}$ the dual space 
of $\frak{g}$ and by $G$ the connected and simply connected Lie group with Lie algebra $\frak{g}$.
Let
$\frak{g}_0\subset \frak{g}_1\subset\ldots\subset \frak{g}_m=\frak{g}$
be  a flag of $\frak{g}$ 
($dim\, {\frak{g}}_i=i$) such that $[\frak{g},\frak{g}_i]\subseteq \frak{g}_{i-1}$ for all $i$
in $\{1,\ldots,m\}$. 
Let also $B=(X_1,...,X_m)$ be a Jordan-H\"older basis adapted to ${(\frak{g}_i)}_i$ that is 
$\frak{g}_i={\R}X_1\oplus\ldots\oplus{\R}X_i$ for all $i$.
Then,
$\frak{g}^{\ast}$ (or $\O$) has a natural layering which can be summarized
 as follows (see also [ACG, Bon, Puk, Ver]).\par\noindent
For $\mu$ in $\frak{g}^\ast$, we define the set of indexes $J_{\mu}=\{j:X_j\notin \frak{g}_{j-1}+\frak{g}_{\mu}\}$,
where $\frak{g}_{\mu}=\{X\in \frak{g}:\forall Y\in \frak{g}, <\mu,[X,Y]>=0\}$.
If $J_{\mu}=\{j_1<\ldots <j_{2r}\}$, we shall have
$\frak{g}=\frak{g}_\mu\oplus {\R}X_{j_1}\ldots\oplus{\R} X_{j_{2r}}$.
Let $\D=\{J_\mu,\mu\in \frak{g}^{\ast}\}$. For $e$ in $\D$, we define the following layer
$$\O_{B}^{e}=\{\mu\in \frak{g}^{\ast}:J_\mu=e\}.$$
By construction, each layer is a $G$-invariant subset of $\frak{g}^\ast$ and
$\frak{g}^\ast$ (resp. $\O$) is a disjoint finite union of layers 
$\frak{g}^\ast=\union_{e\in \D} \O_B^{e}$ (resp. $\O=\union_{e\in \D} (\O_{B}^{e}\cap \O)$).
Note that all the orbits contained in a given layer have 
the same dimension ($card \,e)$.
\par
Now, let $\O_{B}^{e}$ be an arbitrary layer of $\frak{g}^{\ast}$ and 
assume that the orbits contained in $\O_B^{e}$ are $2r$-dimensional.
Then, there exists on $\frak{g}^{\ast}$\par
(i) $m-2r$ polynomial functions $z_1,...,z_{m-2r}$\par
(ii) $2r$ rational
functions $p_1,...,p_r,q_1,...,q_r$ which are regular on $\O_{B}^{e}$
\par\noindent
such that:
\par\noindent
- the polynomial functions $z_1,..., z_{m-2r}$ separate the orbits contained in $\O_B^{e}$;
\par\noindent
- for each orbit $O$ contained in $\O_{B}^e$, there is a diffeomorphism 
(a global Darboux chart)
$\v:O\rightarrow {\R}^{2r}$ of symplectic manifolds defined by the functions $p_i,q_j$.\par
The first layer, noted $V_B$, has additional properties:
it is a Zariski dense open subset of $\frak{g}^\ast$, it contains only orbits of maximal dimension $2d$ and 
the polynomial functions $z_1,..., z_{m-2d}$ separating the orbits of $V_B$ are $G$-invariant. 
Moreover, if we identify the symmetric algebra $S(\frak{g})$ of $\frak{g}$ with the space of polynomial functions on 
$\frak{g}^{\ast}$ and denote by $S(\frak{g})^{G}$ the subring of $S(\frak{g})$ of the $G$-invariant polynomial
functions, then the quotient field of $S(\frak{g})^{G}$ coincides exactly with the field 
$\R(z_1,..., z_{m-2d})$ of rational functions in the $z_k$ variables.
The open set $V_B$ is usually called the generic  set associated to the basis $B$,
the orbits contained in $V_B$ are the generic orbits and 
the corresponding polynomial functions $z_1,...,z_{m-2d}$
are the generic invariants. Since the symplectic foliation of $V_{B}$ is a product foliation 
whose leaves are contractible, it follows immediately from Theorem 3 that the
TP-cohomology of $V_B$ is trivial in degree superior to zero.
The next proposition claims it is even possible to get rid of the choice of the basis $B$.
\begin{Pro}
Let $\frak{g}$ be a nilpotent Lie algebra.
Consider the union $\dis\union_{B} V_B$ 
of all possible generic sets associated to Jordan-H\"older bases $B$ of $\frak{g}$
and denote by $\L$ its regular Poisson structure.
 Then,
$$H^0_{\L,tan}(\dis\union_{B} V_B)=I(\dis\union_{B} V_B)
\quad \hbox{and}\quad H^q_{\L,tan}(\dis\union_{B} V_B)=\{0\}\quad \forall q>0.$$
\end{Pro}
\noindent\par
{\sl Proof.}
We observe first that
$\Pi:{\dis\union_{B} V_B}\rightarrow ({\dis\union_{B} V_B})/ G$ is a locally trivial fibration thus a submersion.
Now, since the orbits contained in $\dis\union_{B} V_B$ 
are connected and cohomologically trivial, the result directly comes from 
Theorem 2 (Section 2).\smallbreak

\begin{Rem}
It is known that $V_B$ is in general strictly included in the set $\O$ of all maximal dimensional coadjoint
orbits (see [SG]). Unfortunately, 
$\dis\union_{B} V_B$ can also be strictly smaller than  $\O$.
For instance, in the case of the filiform Lie algebras (defined in [CG] or [GK]), 
all the $V_B$ coincide and are distinct from $\O$.
\end{Rem}

Now, the following result is very convenient and quite efficient for many examples.
\begin{Pro} 
Let $\frak{g}$ be an m-dimensional nilpotent Lie algebra, $\frak{g}^{\ast}$ the dual space of $\frak{g}$ and $G$ 
the connected and
simply connected Lie group with Lie algebra $\frak{g}$. 
Denote by $\O$ the union of all the coadjoint orbits of maximal dimension $(2d)$.
Still denote by $B=(X_i)$  a Jordan-H\"older basis of $\frak{g}$,
by $V_{B}$  the generic set associated to the basis $B$ and by
$z_1,...,z_{m-2d}$ the generic invariants separating the orbits of $V_B$.
If  $\tilde{\O}$ is any open subset of $\O$ such that the polynomial
 functions $z_1,..., z_{m-2d}$ separate the orbits contained in $\tilde{\O}$  and that the vectors 
$dz_1(\mu),\ldots,dz_{m-2d}(\mu)$ are linearly independent for all $\mu$ in $\tilde{\O}$,
then the TP-cohomology of 
$\tilde{\O}$ (endowed with its linear Poisson structure $\L$) is given by
$$H^0_{\L,tan}(\tilde{\O})=I(\tilde{\O})\quad \hbox{and}\quad H^q_{\L,tan}(\tilde{\O})=\{0\}\quad \forall q>0.$$
\end{Pro}\noindent\par
{\sl Proof.}
 Let us consider the smooth mapping
\begin{align*}
f:&\tilde{\O}\times \tilde{\O} \longrightarrow {\R}^{m-2d}\cr
&(\mu,\eta) \longmapsto (z_1(\mu)-z_1(\eta),\ldots,z_{m-2d}(\mu)-z_{m-2d}(\eta)).
\end{align*}\par\noindent
Note that $\Sigma=f^{-1}(0)$ is not empty (it contains the diagonal set
$\D=\{(\mu,\mu):\mu\in \tilde{\O}\}$).
Moreover,
for all $(\mu,\eta)$ in $\Sigma$, the rank of the linear mapping $f_{{\ast}(\mu,\eta)}$ is, by assumption,
equal to $m-2d$. As a result, $\Sigma $ is a closed submanifold of 
$\tilde{\O}\times \tilde{\O}$. This exactly means (see [Die] p.58) that the space of leaves $\tilde{\O}/ G$
is Hausdorff and that the canonical projection 
$\Pi: \tilde{\O} \rightarrow \tilde{\O}/ G$ is a submersion. The result is
thus again a straightforward consequence of Theorem 2.

\begin{Rem}
If some $\tilde{\O}$ satisfies the conditions of Proposition 3 for a particular Jordan-H\"older basis, it
satisfies these conditions for any Jordan-H\"older basis.
Indeed, since the generic invariants $z_1,..., z_{m-2d}$ associated to any basis $B$ generate the quotient
field of 
$S(\frak g)^G$, the fact that the coadjoint orbits contained in $\tilde{\O}$ are separated or not by 
the polynomial functions $z_1,..., z_{m-2d}$, and 
therefore Proposition 3, do not depend on the choice of the basis $B$.
\end{Rem}

As one can see by studying the examples of Pedersen [Pe1, Pe2], there are some nilpotent Lie algebras for which
Proposition 3, and more generally Theorem 2, cannot be applied to compute the TP-cohomology of the union $\O$
of all maximal dimensional coadjoint orbits.
To deal with these cases which are actually the most fascinating, we propose first to examine the example of
$\frak{g}=\frak{g}_{4,1}$.  
The brackets of this filiform Lie algebra are
$$[X_4,X_3]=X_2,\quad [X_4,X_2]=X_1.$$
Let us identify $\frak{g}^\ast$ with ${\R}^4$ by means of 
the coordinates system $(x_i)$ of $\frak{g}^\ast$ associated to the basis $(X_i)$. 
The 2-dimensional orbits in $\frak{g}^\ast$ are of two kinds.
There are first the orbits of the points $\mu=(\mu_1, \mu_2, \mu_3, \mu_4)$ with $\mu_1\not=0$, which are
parabolic cylinders of the form
$$O_{\mu}=\{(\mu_1, s, \frac{2 \mu_1\mu_3-\mu_2^2+s^2}{2 \mu_1}, \mu_1 t): (s,t)\in {\R}^2\}.$$
Moreover, for the limiting case ($\mu_1=0$),
there are the orbits of the points $\mu=(0, \mu_2, \mu_3, \mu_4)$ with $\mu_2\not=0$, which are 
affine varieties of the form
$$O_{\mu}=\{(0, \mu_2, s, \mu_2 t):(s,t)\in {\R}^2\}.$$
In this example, the regular Poisson manifold $(\O,\L)$ is thus the set 
$$\O=\{(x_1, x_2, x_3, x_4)\in \frak{g}^\ast: x_1^2+x_2^2\not=0\},$$
endowed with the regular Poisson structure $\L$ coming from the Lie bracket.
It is  clear that the generic invariants associated to the basis $(X_i)$, namely $z_1=x_1$ and 
$z_2=x_2^2-2 x_1 x_3$, do not separate the 2-dimensional orbits. 
In fact, the space of leaves $\O/G$ (where $G$ stands for the Lie group corresponding to $\frak{g}$)
equipped with the quotient topology is not Hausdorff. To see this, 
consider the orbits 
$O^{+}$ and $O^{-}$ of the points $(0,1,0,0)$ and $(0,-1,0,0)$.
Denote by $\Pi$ the projection
$\Pi: \O\rightarrow \O/G$.
Let $U$ (resp. $V$) be any neighborhood of $O^{+}$ (resp. $O^{-}$).
Then 
$\Pi^{-1}(U)$ (resp. $\Pi^{-1}(V)$) is an open subset of $\O$ cointaining $(0,1,0,0)$ (resp. $(0,-1,0,0)$). 
Moreover,
$\Pi^{-1}(U)\cap \Pi^{-1}(V)$ intersects the orbits of the points $(\a,1,0,0)$ for sufficiently
small $\a$. It follows that $U\cap V=\Pi(\Pi^{-1}(U)\cap \Pi^{-1}(V))$ is not empty.
\vskip0,5cm

Let us now study the TP-cohomology of $(\O,\L)$.
As always, $$H^0_{\L,tan}(\O)=I(\O).$$
To describe $H_{\L,tan}^1(\O)$, we shall observe that any 
tangential vector field $X$  can be written in the form
$X=a \p_4+ b H_{x_4}$ where $a,b$ are in $C^{\infty}(\O)$
and  that, for such a $X$,  $\s(X)=0$ if and only if 
$$(*)\quad H_{x_4}(b)+\p_4(a)=x_2\p_3(b)+x_1\p_2(b)+\p_4(a)= 0.$$
Let us reduce the study of $H^1_{\L,tan}(\O)$ to the resolution of the partial differential equation
$$H_{x_4}(g)=a_0$$
when $a_0$ and $g$ depend only on the variables $x_1,x_2,x_3$.\par\noindent
Suppose that $X=a \p_4+b H_{x_4}$ satisfies $\s(X)=0$, that is $(\ast)$ holds.
Then, there exists $f$ in $C^{\infty}(\O)$ such that $X=\s f$ if and only if
$$X(dx_i)=\s f(dx_i)=-\{f,x_i\}\quad\forall i,$$
or equivalently, if and only if
$$b=-\p_4(f)\quad\hbox{and}\quad  a=H_{x_4}(f).$$
Thus, the existence of $f$ is equivalent to the existence of $g(x)=g(x_1,x_2,x_3)$ such that
$$f(x)=-\int_0^{x_4} b(x_1,x_2,x_3,t) dt + g(x_1,x_2,x_3)$$
and 
$$H_{x_4}(f)(x)=-\int_0^{x_4} H_{x_4}b(x_1,x_2,x_3,t) dt+ H_{x_4}(g)(x)=a(x).$$
Because of $(*)$, it exactly means that
\begin{align*}
H_{x_4}(f)(x)&=\int_0^{x_4} \p_4(a)(x_1,x_2,x_3,t) dt + H_{x_4}(g)(x) \cr
& =a(x_1,x_2,x_3,x_4)-a(x_1,x_2,x_3,0)+ H_{x_4}(g)(x) \cr
& =a(x_1,x_2,x_3,x_4),\end{align*}
and thus, as announced, that $H_{x_4}(g)(x)=a(x_1,x_2,x_3,0)$.
We want now to prove

\begin{Lem}
Let $a$ be a function in $C^{\infty}(\O)$ depending only on the variables $x_1,x_2,x_3$.
Assume there exists a function $g$ in $C^{\infty}(\O)$, depending also on the variables $x_1,x_2,x_3$,
such that 
$$H_{x_4}(g)=x_2\p_3(g)+x_1\p_2(g)=a.$$
Then, 
$$\lim_{x_1\rightarrow 0}\int_{-1}^1 a(x_1,s, -\frac{1}{2 x_1}+ \frac{s^2}{2 x_1}) 
\frac{ds}{x_1}$$
exists.
\end{Lem}\noindent\par
{\sl Proof.}
Using the change of variables on the open set $U=\{x\in \O: x_1\not=0\}$,
$$u=x_3-\frac{x_2^2}{2 x_1}, v=\frac{x_2}{x_1}, w=x_1,$$
we can see that $g$ is necessarily of the form
\begin{align*}
g(x_1,x_2,x_3)&=\int_0^{\frac{x_2}{x_1}} a(x_1,x_1 t,x_3-\frac{x_2^2}{2x_1}+\frac{t^2 x_1}{2})
dt+c(x_1,x_3-\frac{x_2^2}{2x_1})\cr
 &=\int_0^{x_2} a(x_1,s,x_3-\frac{x_2^2}{2x_1}+\frac{s^2}{2 x_1})
\frac{ds}{x_1} +c(x_1,x_3-\frac{x_2^2}{2x_1}).\end{align*}
\noindent
Thus,
\begin{align*}
g(x_1,1,0)=&\dis\int_0^1 a(x_1,s,-\frac{1}{2 x_1}+\frac{s^2}{2 x_1})\frac{ds}{x_1} + 
c(x_1,-\frac{1}{2 x_1})\cr
\hbox{and}\quad g(x_1,-1,0)=&\dis\int_0^{-1} a(x_1,s,-\frac{1}{2 x_1} +\frac{s^2}{2x_1}) 
\frac{ds}{x_1} + c(x_1,-\frac{1}{2 x_1}).
\end{align*}
\par\noindent
Therefore, $g(x_1,1,0)-g(x_1,-1,0)=\dis\int_{-1}^{1} a(x_1,s,-\frac{1}{2 x_1}+ \frac{s^2}{2 x_1})
\frac{ds}{x_1}.$ We get the result from the fact that  $g$ is  continuous at the points $(0,1,0)$ and
$(0,-1,0)$.\medbreak
While the leaves of $\O$ are cohomologically trivial, the first TP-cohomology space of $\O$ is very
large. Indeed, we have

\begin{Pro}
Let $\O$ be the set of maximal dimensional orbits associated to $\frak{g}=\frak{g}_{4,1}$
and  $\L$ its regular Poisson structure.
Let $G$ be the connected and simply connected Lie group with Lie algebra $\frak{g}$. Denote by
${S(\frak{g})}^{G}$ (resp. ${\mathcal A}(\frak{g})^G$) the ring of G-invariant polynomial 
(resp. analytic) functions over $\frak{g}^\ast$.
For each $\a>\frac{1}{2}$, denote by $T_{\a}$ and $K_{\a}$ the vector fields defined by 
$T_{\a}=t_{\a} \p_4$ and  $K_{\a}=k_{\a} \p_4$ where
$$t_{\a}=\frac{1}{(x_1^2+x_2^2)^{\a}}\quad \hbox{and}\quad 
k_{\a}=\frac{x_1 \exp(\frac{\a}{x_1^2+x_2^2})}{(x_1^2+x_2^2)^2}.$$
Then, the following assertions hold. \par\noindent
(i) The classes $[T_{\a}]$ generate an infinite dimensional space over ${\R}$;\par\noindent 
(ii) The classes $[K_{\a}]$ are linearly independent not only over ${\R}$ but also
over $S(\frak{g})^{G}$ or over ${\mathcal A}(\frak{g})^G$, {\sl i.e.}
 $H_{\L,tan}^1(\O)$, as a $S(\frak{g})^{G}$-module
or as a ${\mathcal A}(\frak{g})^G$-module, is not finitely generated.
\end{Pro}

\begin{Rem} (and convention)
We use in (ii) the integral domains $S(\frak{g})^{G}$ and ${\mathcal A}(\frak{g})^G$, but, from
a differential geometry point of view, it would be more interesting to
consider the whole ring $I(\O)$ of $G$-invariant smooth functions over $\O$.
In fact, due to the complexity of the non-integral domain $I(\O)$, we do not know the $I(\O)$-module structure of 
$H_{\L,tan}^1(\O)$. Nevertheless we conjecture that it is still not finitely generated. 
 \par\noindent
In the sequel, when we say for some Lie algebra $\frak{g}$ that
a TP-cohomology (or a P-cohomology) space of $\O$
is infinite dimensional, we will mean both infinite dimensional as a vector space 
and not finitely generated as a module over
$S(\frak{g})^{G}$ or ${\mathcal A}(\frak{g})^G$.
\end{Rem}

\noindent\par
{\sl Proof of Proposition 4.}
For all $0<x_1<1$, we have
 $$\dis\int_0^1 \frac{ds}{(x_1^2+s^2)^{\a} x_1}\geq \int_0^{x_1} \frac{ds}{(x_1^2+s^2)^{\a} x_1}
\geq \frac{1}{2^{\a} x_1^{2 \a}}\quad (\forall \a) $$
and 
\begin{align*}
\int_0^1\frac{ ds}{(x_1^2+s^2)^{\a} x_1}&\leq \int_0^{x_1} \frac{ds}{(x_1^2+s^2)^{\a} x_1}
+\int_{x_1}^1\frac{ds}{s^{2\a} x_1}
\leq c_{\a} \frac{1}{x_1^{2\a}}\quad (\forall \a>\frac{1}{2}) 
\end{align*}
where $c_{\a}=\frac{2\a}{2\a-1}.$
Now, assume that for some $p$, $\dis\sum_{i=1}^p \l_{\a_i} [T_{\a_i}]=0$, the
$\l_{{\a}_i}$ being in ${\R}$ and $\a_1<\ldots<\a_p$.
Using Lemma 3, we directly see that 
$$E=\int_0^1\dis\sum_{i=1}^{p} \l_{\a_i} \frac{ds}{(x_1^2+s^2)^{\a_i} x_1}$$
 must have a limit when $x_1$ tends to zero.
But, when $0<x_1<1$,
\begin{align*}
|E|&\geq |\dis\int_0^1\frac{\l_{\a_p} ds}{x_1(x_1^2+s^2)^{\a_p}}|-|\dis\int_0^1\frac{\l_{{\a}_{p-1}} ds}
{ x_1(x_1^2+s^2)^{\a_{p-1}}}|\cr
& \ldots -|\dis\int_0^1\frac{\l_{\a_1} ds}{x_1(x_1^2+s^2)^{\a_1}}|\cr
&\dis\geq \frac{|\l_{\a_p}|}{2^{\a_p} x_1^{2\a_p}}-\frac{|\l_{\a_{p-1}}|c_{\a_{p-1}}}{x_1^{2
\a_{p-1}}}-\ldots -\dis\frac{|\l_{\a_1}| c_{\a_1}}{x_1^{2\a_1}}.
\end{align*}
Thus, $\l_{\a_p}$ must be $0$
and a step-by-step application of the same argument shows that $\l_{\a_i}=0$ for all $i$.
It implies that the classes $[T_{\a}]$
generate an infinite dimensional vector space over ${\R}$. That ends the proof of (i).\par\noindent
To prove (ii), it is enough to see that
\begin{align*}
\int_0^1 \frac{2 \exp(\frac{\a}{x_1^2+s^2})}{(x_1^2+s^2)^2} ds&\geq \dis\int_0^1 
\frac{2s \exp(\frac{\a}{x_1^2+s^2})}{(x_1^2+s^2)^2}ds\cr
&\geq \frac{1}{\a} (\exp(\frac{\a}{x_1^2})-\exp(\frac{\a}{x_1^2+1})).
\end{align*}
and that
$$\int_0^1\frac{2 \exp(\frac{\a}{x_1^2+s^2})}{(x_1^2+s^2)^2} ds\leq
\frac{2\exp(\frac{\a}{x_1^2})}{x_1^4}.$$ Thus we can use the same argument as in (i).
\medbreak
To finish the discussion about this example, we shall prove the vanishing of 
$H^2_{\L,tan}(\O)$.
Let $A$ be a tangential 2-tensor field. Necessarily, 
 $$A=\v \L=\v (x_2\p_4\w\p_3+ x_1 \p_4\w\p_2)$$ for some function $\v$ in $C^{\infty}(\O)$ and  
$\s(A)=0.$ We have thus to find a tangential vector field $B$  such that
$A=\s (B)$, or equivalently, to find $a,b$ in $C^{\infty}(\O)$ such that
$H_{x_4}(b)+\p_4(a)=\v$. We immediately check that 
$$a=\int_0^{x_4}\v(x_1,x_2,x_3,t) dt\quad\hbox{and}\quad b=0$$ 
are convenient. 
Therefore, $H^2_{\L,tan}(\O)=\{0\}$.\par\noindent
This fact can also be
deduced from Proposition 1 (Section 2), which identifies the TP-cohomology with a $\check{C}$ech
cohomology. Indeed, it is easy to see that $\O$ admits a good covering ${\mathcal U}=(U_i)$ 
without any non-trivial intersections of three open sets $U_i$.\par
Now, since the rank of $\O$ is 2, we have $H^k_{\L,tan}(\O)=\{0\}$ for all $k>2$.

\begin{Rem}
An immediate consequence of the above analysis concerning the nilpotent case is that
the tangential star products on $\dis\union_{B} V_B$ (and on any regular Poisson manifold
$\tilde{\O}$ satisfying the conditions of Proposition 3) are all 
 equivalent because of the vanishing of the second TP-cohomology space
(Propositions 2 and 3).
The same is true for the union $\O$ of all maximal dimensional coadjoint orbits in $\frak{g}_{4,1}^{\ast}$. \end{Rem}

\medbreak
The computation of the TP-cohomology spaces for $\frak{g}_{5,4}$,
 $\frak{g}_{6,18}$ or for all the filiform Lie algebras, studied for instance in [Pe1, Pe2, CG, GK, BLM],
leads to the same results as in the $\frak{g}_{4,1}$-case. Due to these examples, we 
believe  that the TP-cohomology spaces of a regular Poisson 
manifold $M$ are huge and rather complicated to compute
whenever the quotient space of $M$ by the foliation is not  
Hausdorff.
To confirm this observation, we shall now study more varied examples. 

\

\section{Further examples}\label{sec4}

Let $(M,\L)$ be a 3-dimensional regular Poisson manifold. If we exclude the trivial
case where $\L=0$, we can suppose $M$ to be of rank 2.
To describe the P-cohomology of $M$, we are going to use Vaisman's notations and computations ([Va2] p.69).
Of course, 
$$H{}^0_{\L}(M)=H{}^{0}_{\L,tan}(M)=I(M).$$
Now, we have by definition

$$\dis H{}^1_{\L}(M)=\frac{\{Q\in {\mathcal V}^1(M): \s(Q)=0\}}
                 {\{\s f:f\in C^{\infty}(M)\}}.$$
Let us then choose a transversal distribution $\nu{\mathcal F}$. Let us also use the decompositions
$\dis\somme_{q}{\mathcal V}_{0,q}(M)$, $\dis\somme_{q}\O_{0,q}(M)$, $\s=\s'+\s''$ and $d=d'+d''+d_{2,-1}$.
Each element  $Q$ of ${\mathcal V}^1(M)$ can thus be written in the form
$$Q=Q_{0,1}+Q_{1,0}$$ where $Q_{0,1}$ and  $Q_{1,0}$ are of type  $(0,1)$ and $(1,0)$ respectively.
Furthermore,  $\s(Q)=0$ if and only if $\s(Q_{0,1}+Q_{1,0})=0$,
or else, if and only if $$\s'(Q_{1,0})+\s''(Q_{1,0})+\s''(Q_{0,1})=0,$$
that is, if and only if $$\s''(Q_{1,0})=0\quad\hbox{and}\quad\s'(Q_{1,0})+\s''(Q_{0,1})=0.$$
Consider now the linear mapping $p$ defined by
\begin{align*}
&p:H{}_{\L}^1(M)\longrightarrow \tilde{{\mathcal V}}_{1,0}(M)\cr
&[A]\longmapsto A_{1,0}
\end{align*}
where $$\tilde{{\mathcal V}}_{1,0}(M)=\{A\in {\mathcal V}_{1,0}(M):\s''(A)=0\}.$$
It follows that
$$H{}_{\L}^1(M)\cong Ker (p)\oplus Im(p).$$
Now \begin{align*}
Ker(p)&=\{[A]\in H{}_{\L}^1(M): A_{1,0}=0\}\cr
 &=\dis\frac{{\{A\in {\mathcal V}_{0,1}(M)} : \s(A)=0\}}
{\{\s f:f\in C^{\infty}(M)\}}\cr
&=\dis\frac{{\{A\in{\mathcal V}_{0,1}(M)}:\s''(A)=0\}}
{\{\s'' f:f\in C^{\infty}(M)\}}\cr
&=H{}_{\L,tan}^{1}(M).
\end{align*}
Moreover,
$$Im(p)=\{A\in {\mathcal V}_{1,0}(M):\s''(A)=0 \quad\hbox{and}\quad
\exists B \in {\mathcal V}_{0,1}(M)/ \s(A)+ \s''(B)=0 \}.$$
Let us compute the second order space:
\begin{align*}
H{}_{\L}^2(M)&=\dis\frac{\{Q\in {\mathcal V}^2(M):\s(Q)=0\}}
{\{\s(V):V\in{\mathcal V}^1(M)\}}\cr
&=\dis\frac{\{Q=Q_{0,2}+Q_{1,1}\in{\mathcal V}^2(M): \s''(Q_{1,1})=0\}}
{\{\s(V): V=V_{0,1}+V_{1,0}\in{\mathcal V}^1(M)\}}\cr
&=\dis\frac{\{ Q_{0,2}\in {\mathcal V}_{0,2}(M)\}}
 {\{\s''(V_{0,1})+\s'(V_{1,0})\}}\oplus
\dis\frac{\{Q_{1,1}\in{\mathcal V}_{1,1}(M):\s''(Q_{1,1})=0\}}
{\{\s''(V_{1,0})\}}\cr
&=\dis\frac{\dis\frac{\{Q_{0,2}\}}
 {\{\s''(V_{0,1})\}}}{\s({\mathcal V}_{1,0}(M))}\oplus
\dis\frac{\{Q_{1,1}:\s''(Q_{1,1})=0\}}{\{\s''(V_{1,0})\}}\cr
&=\frac{H^2_{\L,tan}(M)}{\s({\mathcal V}_{1,0}(M))} \oplus 
\frac{ Ker (\s''_{| {\mathcal V}_{1,1}(M)})}{\s''({\mathcal V}_{1,0}(M))}.
\end{align*}
In the same way, we get the third order space:
$$ H{}_{\L}^3(M)=\dis\frac{\{Q_{1,2}\}}{\{\s(V_{0,2}+V_{1,1})\}}
=\dis\frac{\{Q_{1,2}\}}{\{\s''(V_{1,1})\}}=\frac{ {\mathcal V}_{1,2}(M)}{\s''({\mathcal
V}_{1,1}(M))}.$$

\vskip0,5cm
 As seen in Section 2, for all regular Poisson manifold $(M,\L)$, 
$(\dis\somme_{q}\O_{0,q}(M),d'')$ and $(\dis\somme_{q}{\mathcal V}_{0,q}(M),\s'')$
are isomorphic as complexes of $I(M)$-modules.
In fact, it is always possible to define an isomorphism between 
$\O_{p,q}(M)$ and ${\mathcal V}_{p,q}(M)$ for all $p$ and $q$.
For this, one can consider, as Vaisman did in [Va2], an Euclidean metric on $\nu^{\ast}{\mathcal F}$.
That leads to an isomorphism between $\nu^{\ast}{\mathcal F}\oplus T^{\ast}{\mathcal F}$
and $\nu{\mathcal F}\oplus T{\mathcal F}$, which can  naturally be extended to the required isomorphism
between $\O_{p,q}(M)$ and ${\mathcal V}_{p,q}(M)$. In general, however, this isomorphism is not an
isomorphism of complexes. Via the next proposition, we give a simple situation 
where  $(\dis\somme_{q}\O_{p,q}(M),d'')$ and 
$(\dis\somme_{q}{\mathcal V}_{p,q}(M),\s'')$ are effectively isomorphic complexes for all $p$.
\begin{Pro}
Let $(M,\L)$ be a regular Poisson manifold. Denote by ${\mathcal F}$ the symplectic foliation of $M$.
Choose a transversal distribution $\nu{\mathcal F}$ and
assume that $\O_{1,0}(M)$ and ${\mathcal V}_{1,0}(M)$ are isomorphic free $I(M)$-modules with bases $(\b_i)$ and $(X_i)$
respectively  ($\b_i$ and $X_i$ being globally defined).  
Suppose also that $d''\b_i=0$ and  $\s''(X_i)=0$.
Then, for all $p$, $(\dis\somme_{q}\O_{p,q}(M),d'')$ and $(\dis\somme_{q}{\mathcal V}_{p,q}(M),\s'')$
are isomorphic complexes of $I(M)$-modules. In particular, $H^q(M,\P^p({\mathcal F}))$ and 
$H^q(\dis\somme_{k}{\mathcal V}_{p,k}(M),\s'')$ are isomorphic as $I(M)$-modules.
\end{Pro}\noindent\par
{\sl Proof.}
Recall that $\#:T^{\ast}M\rightarrow TM$ can be extended to an $I(M)$-modules isomorphism
$\tilde{\#}:\O_{0,q}(M)\rightarrow {\mathcal V}_{0,q}(M)$ satisfying $\s ''\circ \tilde{\#}=-\tilde{\#}\circ d''$.
Consider now the mapping 
$\hat{\#}:\O_{p,q}(M)\rightarrow {\mathcal V}_{p,q}(M)$ defined
by 
$$\hat{\#}(\sum_{i_1,...,i_p} \a_{i_1,...,i_p}\w\b_{i_1}\w\ldots \w\b_{i_p})=
\dis\sum_{i_1,...,i_p} (\tilde{\#}(\a_{i_1,...,i_p})\w X_{i_1}\w\ldots \w X_{i_p})$$
where the $\a_{i_1,...,i_p}$ are in $\O_{0,q}(M)$.
It is not difficult to check that $\hat{\#}$ is an $I(M)$-modules isomorphism 
and that $\s ''\circ \hat{\#}=-\hat{\#}\circ d''$. This ends the proof.

\medbreak
The following definitions are quite standard, we recall them for completeness.
Let $V\rightarrow B$ be a vector bundle whose fibers are $q$-dimensional. We shall say that
$V\rightarrow B$ (or simply $V$) is orientable if
the bundle $\L^q V\rightarrow B$ admits a global nonsingular ({\sl i.e.} nowhere vanishing) section.
If $V\rightarrow B$ is orientable, so is its dual $V^\ast\rightarrow B$. Recall also that a manifold
$M$ is said to be orientable if $TM$ (or $T^{\ast}M$) is orientable.\par
Moreover, we shall say that a foliation ${\mathcal F}$ on $M$ is (co)orientable if its normal
bundle $N{\mathcal F}$ is orientable. 
In the important case of a 1-codimensional foliation ${\mathcal F}$, this foliation ${\mathcal F}$ is orientable if
and only if there exists a nonsingular 1-form $\b$  vanishing
exactly on vectors tangent to the leaves of ${\mathcal F}$.
In this case, we say that $\b$ defines the foliation.\par
Note that neither the leaves nor the total manifold $M$ of
an orientable foliation need to be orientable.
However, if ($M,\L$) is a 
regular Poisson manifold (of rank $2n$) and ${\mathcal F}$ is the symplectic foliation of $M$, 
then the situation is somewhat simpler.
Since the tensor $\L^n$ defines a global nonsingular section of $\L^{2n}T{\mathcal F}$, the tangent bundle $T{\mathcal F}$
of ${\mathcal F}$ is orientable. In other words, the symplectic foliation of a regular Poisson manifold $M$ 
is orientable if and only
if $M$ is orientable. \par

\medbreak  
We turn back now to the case where $(M,\L)$ is 3-dimensional and give a result  which will 
be useful later.\par\noindent

\begin{Pro}
Let $(M,\L)$ be a regular Poisson manifold of dimension 3 and rank 2.
Denote by ${\mathcal F}$ the symplectic foliation of $M$. 
Suppose that $M$ is orientable and that one can choose the 1-form $\b$ defining ${\mathcal F}$ 
 such that $d\b=0$. Then
\par\noindent
(i) the second and third $P$-cohomology spaces of $M$ are
 \begin{align*}
H_{\L}^2(M)&\simeq \frac{ H^2_{\L,tan}(M)}{\s({\mathcal V}_{1,0}(M))} \oplus H^1(M,\P^1({\mathcal
F}))\cr
 H_{\L}^3(M)&\simeq H^2(M,\P^1({\mathcal F})).\end{align*}
(ii) for all $q$, $H^q_{\L,tan}(M)$ and 
$H^q(M,\P^1({\mathcal F}))$ are isomorphic as $I(M)$-modules.
\end{Pro}
\noindent\par
{\sl Proof.}
Let us identify the normal bundle $N{\mathcal F}$ of ${\mathcal F}$ with a sub-bundle $\nu{\mathcal F}$
of $TM$. Since $M$ (or ${\mathcal F}$) is orientable, there exists a nonsingular vector
field $X$ (globally defined) such that, for each $x$, 
$\nu_{x}{\mathcal F}$ (resp. $\nu_{x}^{\ast}{\mathcal F}$)
is spanned by $X_x$ (resp. $\b_x$).
By construction,  $\s''X=0$ and $d\b=d''\b=0$. Moreover,
 $\O_{1,0}(M)$ and ${\mathcal V}_{1,0}(M)$ are isomorphic free and with basis
$(\b)$ and $(X)$ respectively.
The point (i) is thus a direct corollary of Proposition 5.\par\noindent
To prove (ii), denote by $\Phi:H^q(M,\P^0({\mathcal F}))\rightarrow H^q(M,\P^1({\mathcal F}))$ the mapping
defined by $\Phi([\a])=[\a\w\b]$ for all $\a$ in $\O_{0,q}(M)$ such that $d''\a=0$.
Clearly, $\Phi$ is both well defined and  bijective.
Therefore, $H^q(M,\P^0({\mathcal F}))$ (which is isomorphic to
$H^q_{\L,tan}(M)$ by Theorem 1) coincides with $H^q(M,\P^1({\mathcal F}))$.
This ends the proof.

\medbreak

As we already said, each Lie algebra provides a natural regular Poisson manifold:
the union $\O$ of all maximal dimensional coadjoint orbits.
Let us now study the TP-cohomology and the P-cohomology 
of $\O$ for any 3-dimensional Lie algebra.\medbreak
First recall that every non-abelian  3-dimensional Lie algebra  is isomorphic to exactly one in the following 
list: (see [Br] for instance)\par\noindent
- a nilpotent Lie algebra, namely the 3-dimensional Heisenberg algebra, given by: 
$[X_1,X_2]=X_3$;\par\noindent
- a solvable non exponential Lie algebra, namely $\frak{e}(2)$, defined by:\par\noindent
$[X_1,X_2]=-X_3,\, [X_1,X_3]=X_2$;\par\noindent
-several (solvable) exponential Lie algebras:\par\noindent
$\ast$ the algebra spanned by $X_1,X_2,X_3$ with 
$[X_1,X_2]=X_2$; \par\noindent
$\ast$ the ``book algebra":\par\noindent
 $[X_1,X_2]=X_2,\, [X_1,X_3]=X_3$;\par\noindent
$\ast$ ``Gr\'elaud's Lie algebras":\par\noindent
$[X_1,X_2]=X_2-\s X_3, \, [X_1,X_3]=X_3+\s X_2$  where $\s>0$;\par\noindent
$\ast$ $[X_1,X_2]=X_2,\,[X_1,X_3]=\frac{1}{\tau} X_3$ where $\tau>1$;\par\noindent
$\ast$ $[X_1,X_2]=X_2+X_3,\, [X_1,X_3]=X_3$;\par\noindent
$\ast$ $[X_1,X_2]=X_2,\, [X_1,X_3]=-X_3$;\par\noindent
$\ast$ $[X_1,X_2]=X_2,\, [X_1,X_3]=\frac{1}{\tau} X_3$ where $\tau<-1$;\par\noindent
-two simple Lie algebras:\par\noindent
$\frak{su}(2)$ and $\frak{sl}(2)$.\par\noindent
Among them, there are some easy examples. Let us pass them in review.
\medbreak

\subsection{The 3-dimensional Heisenberg algebra}

\

This nilpotent Lie algebra is defined by the bracket $[X_1,X_2]=X_3$.\par\noindent
For this example, the non-trivial orbits are planes and the 
regular Poisson manifold $(\O,\L)$ is
the set \par\noindent
$$\O=\{(x_1,x_2,x_3): x_3\not=0\},$$
endowed with its regular Poisson structure $\L$.
We identify $\O$
with ${\R}^2\times {\R}^{\ast}$ by means of the  Weinstein chart [SG, Wei]
mapping an element 
$x=(x_1,x_2,x_3)$ of $\O$ to $(p=x_1,q=\frac{x_2}{x_3},z=x_3).$ 
\par\noindent
Using Theorem 3 (or Theorem 2), we see that the TP-cohomology of $\O$ is trivial in degree superior to zero.\par\noindent
Furthermore, $\L$ is transversally constant with respect to the transversal distribution $\nu{\mathcal F}=T 
{\R}^{\ast}$.
Using Theorem 4 or by direct computation, we get the P-cohomology of $\O$:
\begin{align*}
H^0_{\L}(\O)&=H_{\L,tan}^0(\O)=I(\O)=C^{\infty}({\R}^{\ast})\cr
H^1_{\L}(\O)&\simeq \{u \, dz: u\in
I(\O)\}\simeq I(\O)\cr
 H^k_{\L}(\O)&= \{0\}\quad\forall k>1.
\end{align*}\medbreak
\subsection{The Lie-algebra $\frak{e}(2)$ of the Euclidean 2-dimensional group}

\

The brackets are:
$$[X_1,X_2]=-X_3,\quad [X_1,X_3]=X_2.$$
The non-trivial orbits are cylinders $C_r$ with radius $r>0$, and
 the regular Poisson manifold $(\O,\L)$ associated to this Lie algebra is the set
$$\O=\{(x_1,x_2,x_3)\in \frak{e}(2)^{\ast}:x_2^2+x_3^2\not=0\},$$
endowed with its regular Poisson structure $\L$.
We identify $\O$ with ${\R}\times {\T}\times {\R}_{+}^{\ast}$ by means of the Weinstein chart
mapping  an element  $x=(x_1,x_2,x_3)$ of $\O$ to the point $(p,q,r)$ defined by
$$p=x_1,\quad 
e^{\imath q }=\frac{x_2+\imath x_3}{\sqrt{x_2^2+x_3^2}},\quad r=\sqrt{x_2^2+x_3^2}.$$
Using Theorem 3, we obtain the TP-cohomology of $\O$:
\begin{align*}
H^0_{\L,tan}(\O)&=I(\O)=C^{\infty}({\R}^{\ast}_{+})\cr
H^1_{\L,tan}(\O)&\simeq H_{DR}^1({\R}\times {\T})\otimes C^{\infty}({\R}^{\ast}_{+})\cr
           &\simeq \{[\rho(r) dq]: \rho(r)\in I(\O)\}\simeq I(\O)\cr
H^k_{\L,tan}(\O)&=\{0\}\quad \forall k>1.
\end{align*}
Now, $\L$ is transversally constant with respect to the transversal distribution $\nu{\mathcal F}=
T{\R}_{+}^{\ast}$.
By Theorem 4, the P-cohomology of $\O$ is:
\begin{align*}
H^0_{\L}(\O)&=I(\O)\cr
H^1_{\L}(\O)&\simeq H^1_{\L,tan}(\O)\oplus \{ u\, dr: u\in I(\O)\}\cr
&\simeq  I(\O)\oplus I(\O)\cr
H_{\L}^2(\O)&\simeq  H_{DR}^1({\R}\times {\T})\otimes \O^1({\R}^{\ast}_{+})\cr
&\simeq \{[\rho(r) dq\w dr]: \rho(r)\in I(\O)\}\simeq I(\O)\cr
H_{\L}^k(\O)&=\{0\}\quad\forall k>2.
\end{align*}\medbreak\par
\subsection{The Lie algebra $\frak{a}$ of the affine group}

\

The only non-vanishing bracket is $[X_1,X_2]=X_2.$\par\noindent
The non-trivial orbits are half planes ($x_3$ and $sign(x_2)$ fixed)
 and the regular Poisson manifold $(\O,\L)$ associated to this Lie algebra is
the set 
$$\O=\{(x_1,x_2,x_3)\in \frak{a}^\ast: x_2\not=0\},$$
endowed with its regular Poisson structure $\L$.
We identify $\O$ with ${\R}^{\ast}\times {\R}\times {\R}$ by means of
the Weinstein chart mapping an element $x=(x_1,x_2,x_3)$ of $\O$ to
the point $(p,q,z)$ defined by
$$p=x_2,\quad q=\frac{-x_1}{x_2},\quad z=x_3.$$
By Theorem 3 (or Theorem 2), the TP-cohomology of $\O$ is trivial in  degree  superior to zero.\par\noindent
Again, $\L$ is transversally constant with respect to the transversal distribution $\nu{\mathcal F}=
T{\R}$.
By Theorem 4, the P-cohomology of $\O$ is:
\begin{align*}
 H^0_{\L}(\O)&=H^0_{\L,tan}(\O)=I(\O)\cr
&\simeq C^{\infty}((\{1\}\times {\R})\cup (\{-1\}\times {\R}))\cr
&\simeq C^{\infty}({\R})\oplus C^{\infty}({\R})\cr
H^1_{\L}(\O)&\simeq\{u\, dz: u\in I(\O)\}
\simeq I(\O)\cr
H_{\L}^{k}(\O)&=\{0\}\quad \forall k>1.
\end{align*}

\medbreak
\subsection{The book algebra}

\

This Lie algebra is given by the following brackets:
$$[X_1,X_2]=X_2,\quad [X_1,X_3]=X_3.$$
The 2-dimensional orbits are characterized by an invariant $\theta$ and are of the form
$$O_{\theta}=\{(s,e^t \cos(\theta), e^t \sin(\theta)): (s,t)\in {\R}^2\}.$$
The corresponding regular Poisson manifold $(\O,\L)$ is thus the set
$$\O=\{(x_1,x_2,x_3): x_2^2+x_3^2\not=0\},$$
with its regular Poisson structure $\L$.
We identify $\O$ with ${\R}\times {\R}\times {\T}$ by means of the Weinstein chart mapping an
element 
$x=(x_1,x_2,x_3)$ 
of $\O$ to the point $(p,q,\theta)$ defined by
$$p=x_1,\quad q=\frac{1}{2} \ln(x_2^2+x_3^2),\quad 
 e^{\imath \theta }=\frac{x_2+\imath x_3}{\sqrt{x_2^2+x_3^2}}. $$
By Theorem 3 (or Theorem 2), the TP-cohomology of $\O$ is trivial in degree superior to zero.\par\noindent
Moreover, $\L$ is transversally constant with respect to $\nu{\mathcal F}=T {\T}$.
By Theorem 4, the  P-cohomology of $\O$ is:
\begin{align*}
H^0_{\L}(\O)&=H_{\L,tan}^0(\O)=I(\O)=C^{\infty}({\T})\cr
H^1_{\L}(\O)&\simeq \{ u\, d\theta:u \in I(\O)\}\simeq I(\O)\cr
H^k_{\L}(\O)&=\{0\} \quad\forall k> 1.
\end{align*}
\medbreak

\

\subsection{Similar to the book algebra: Gr\'elaud's Lie algebras}

\

The brackets of these Lie  algebras, studied by Gr\'elaud in [Gr\'e], are
$$[X_1,X_2]=X_2-\s X_3, \quad [X_1,X_3]=X_3+\s X_2\quad\s>0.$$
The 2-dimensional orbits are a ``spiral" version of those of the book algebra, they are of the form
$$O_{\theta}=\{(s, e^t \cos (\theta +\s t), e^t\sin (\theta+\s t)): (s,t)\in{\R}^2\}$$ 
where  $\theta$ is defined by
$$\dis e^{\imath \theta }=\frac{x_2+\imath x_3}{\sqrt{x_2^2+x_3^2}} e^{-\frac{\imath \s}{2}
\ln(x_2^2+x_3^2)}.$$ The regular Poisson manifold $(\O,\L)$ is thus the set 
$$\O=\{(x_1,x_2,x_3): x_2^2+x_3^2\not=0\},$$
with its regular Poisson structure $\L$.\par\noindent
The situation here is identical to that of the book algebra. First,
 $\O$ can be identified with $\dis{\R}\times {\R}\times {\T}$ and
the TP-cohomology of $\O$ is trivial in degree superior to zero.
Moreover, $\L$ is  transversally constant with respect to 
$\nu{\mathcal F}=T {\T}$
so that the P-cohomology of $\O$ is:
\begin{align*}
H^0_{\L}(\O)&= H^0_{\L,tan}(\O)=I(\O)=C^{\infty}({\T})\cr
H^1_{\L}(\O)&= \{u\, d\theta:u\in I(\O)\}\simeq I(\O)\cr
H^k_{\L}(\O)&=\{0\}\quad\forall k> 1.
\end{align*}\medbreak
\subsection{Other examples very close to the book algebra}

\

Consider the family of Lie algebras defined by the brackets: 
$$[X_1,X_2]=X_2,\quad [X_1,X_3]=\frac{X_3}{\tau }\quad \tau >1.$$
The 2-dimensional orbits can be parameterized with obvious notations by
$$O_{\mu}=\{(s, e^{-t}\mu_{2}, e^{-\frac{t}{\tau }}\mu_{3}): (s,t)\in {\R}^2\},$$
they  are all cohomologically trivial and the regular Poisson manifold $(\O,\L)$ is the set
$$\O=\{(x_1,x_2,x_3): x_2^2+x_3^2\not=0\}$$
with its regular Poisson structure $\L$.\par\noindent
Let us now prove that the symplectic foliation ${\mathcal F}$ of $\O$ is given by a submersion $\Pi$ from $\O$
to the circle ${\bf S}^1$.  
Let $$F: {\R}\times {\R}^2\backslash \{(0,0)\}\rightarrow {\R}$$ be the mapping
defined by 
$$F(t,x_2,x_3)=e^{-2t} x_2^2+e^{-\frac{2t}{\tau}} x_3^2-1.$$
Using the standard Implicit Function Theorem, we see
there is a unique smooth function $\v:{\R}^2\backslash \{(0,0)\}\rightarrow {\R}$ of the variables
$x_2,x_3$ such that 
$$F(\v(x_2,x_3),x_2,x_3)=0\quad\forall (x_2,x_3)\in {\R}^2\backslash \{(0,0)\}.$$
Moreover, the partial derivatives of $\v$ are 
\begin{align*} 
\frac{\p \v}{\p x_2}&=\frac{ e^{-2 \v}x_2}{e^{-2\v}x_2^2+e^\frac{-2\v}{\tau }
\frac{x_3^2}{\tau }}\cr
\frac{\p \v}{\p x_3}&=\frac{ e^{\frac{-2  \v}{\tau }}x_3}{e^{-2\v}x_2^2+e^\frac{-2\v}{\tau
}\frac{x_3^2}
{\tau }} .\end{align*} 
It is  easy to check that the map
$\Pi:\O\rightarrow {\bf S}^1$ defined by
$$(x_1,x_2,x_3)\longmapsto \Pi(x_1,x_2,x_3)=(e^{-\v(x_2,x_3)} x_2, e^\frac{-\v(x_2,x_3)}{\tau} x_3)$$
is a submersion.
Thus, by Theorem 2, the TP-cohomology of $\O$ is trivial in degree superior to zero
just like in the book algebra example.\par\noindent
Furthermore, observe that $\O$ is orientable and that there is 
a 1-form $\b$ defining the foliation ${\mathcal F}$ such that $d\b=0$, namely 
$\b=d\theta$ where $\theta$ is given by
$$e^{\imath \theta} =e^{-\v}x_2+\imath e^\frac{-\v}{\tau}x_3.$$
Therefore, by Proposition 6, the P-cohomology of $\O$
is the same as in the case of the book algebra:
\begin{align*}
H^0_{\L}(\O)&= H^0_{\L,tan}(\O)=I(\O)\cr
H^1_{\L}(\O)&\simeq \{u\, d\theta:u\in I(\O)\}\simeq I(\O)\cr
H^k_{\L}(\O)&\simeq \{0\}\quad\forall k> 1.
\end{align*}
\medbreak
An analogous example is  the Lie algebra defined by the brackets:
$$[X_1,X_2]=X_2+X_3,\quad [X_1,X_3]=X_3.$$
In this case, the 2-dimensional orbits can be parameterized by
$$O_{\mu}=\{(s, e^{-t}(\mu_{2}-\mu_{3}t),e^{-t} \mu_{3}): (s,t)\in {\R}^2\},$$
they are all cohomologically trivial and the regular Poisson manifold $(\O,\L)$ is the
set
$$\O=\{(x_1,x_2,x_3): x_2^2+x_3^2\not=0\}$$
with its regular Poisson structure $\L$.\par\noindent
Now, let $G$ denote the following map
$$G(t,x_2,x_3)=e^{-2t}((x_2-x_3 t)^2+x_3^2)-1.$$
By the Implicit Function Theorem, the equation $G(t,x_2,x_3)=0$ determines $t$ implicitely as a smooth function 
$\psi:{\R}^2\backslash \{(0,0)\}\rightarrow {\R}$ of the variables $x_2, x_3$. The partial 
 derivatives of $\psi$ are
\begin{align*} 
\frac{\p \psi}{\p x_2}&=\frac{ x_2-\psi x_3}{(x_2-(\psi-\frac{1}{2}x_3)^2+ \frac{3}{4} x_3^2}\cr
\frac{\p \psi}{\p x_3}&=\frac{ x_3+\psi^2 x_3-\psi x_2}{(x_2-(\psi-\frac{1}{2} x_3)^2+\frac{3}{4}
x_3^2} .\end{align*} 
As before, the symplectic foliation ${\mathcal F}$ of $\O$ is given by a submersion, namely the map
$\Pi:\O\rightarrow {\bf S}^1$ defined by
$$\Pi(x_1,x_2,x_3)=(e^{-\psi(x_2,x_3)}(x_2-x_3 \psi(x_2,x_3) ), e^{-\psi(x_2,x_3)}x_3).$$ 
Thus, by Theorem 2, the TP-cohomology of $\O$
is trivial
in degree superior to zero.\par\noindent
Moreover, $\O$ is orientable and there exists a 1-form $\b$ defining ${\mathcal F}$ such that
$d\b=0$, namely
$\b=d\theta$ 
where $\theta$ is now defined by
$$ e^{\imath \theta}=e^{-\psi}((x_2-x_3\psi)+\imath x_3).$$
The P-cohomology of $\O$ is thus again:
\begin{align*}
H^0_{\L}(\O)&= H^0_{\L,tan}(\O)=I(\O)\cr
H^1_{\L}(\O)&\simeq \{u\, d\theta:u\in I(\O)\}\cr
&\simeq I(\O)\cr
H^k_{\L}(\O)&\simeq\{0\}\quad\forall k> 1.
\end{align*}

\medbreak

\subsection{The simple Lie algebra $\frak{su}(2)$}

\

This Lie algebra is defined by the following brackets:
$$[X_1,X_2]=X_3,\quad [X_2,X_3]=X_1,\quad [X_3,X_1]=X_2.$$
The non-trivial orbits are 2-spheres and the regular Poisson manifold $(\O,\L)$ is 
$\O=\frak{su}(2)^{\ast} \backslash \{0\}$ endowed with its regular Poisson structure $\L$.
We naturally identify $\O$ with 
${\bf S}^2\times {\R}^{\ast}_{+}$.
Using Theorem 3, we get the TP-cohomology of $\O$: 
\begin{align*}
H_{\L,tan}^0(\O)&=I(\O)=C^{\infty}({\R}^{\ast}_{+})\cr
H_{\L,tan}^1(\O)&=\{0\}\cr
H_{\L,tan}^2(\O)&=\{[u \L]:u\in I(\O)\}\cr
&\simeq\{[u \omega_{\L}]: u\in I(\O)\}\simeq I(\O) 
\end{align*}
where $\omega_{\L}$ denotes the foliated 2-form associated to $\L$ {\sl i.e.}
$$\omega_{\L}=\frac{x_3 dx_1\w dx_2+ x_1 dx_2\w dx_3 + x_2 dx_3\w dx_1}{r^2}\quad
(r=\sqrt{x_1^2+x_2^2+x_3^2}).$$ As for the other TP-cohomology spaces,
$$H^k_{\L,tan}(\O)=\{0\}\quad\forall k>2.$$\noindent
It is important to note that, in this example, $\O$ is a product
 $S\times R$ like in Theorem 4. However, 
$\L$ is not transversally constant so that Theorem 4 cannot be used.\par\noindent 
We shall here compute the P-cohomology of $(\O,\L)$.
First, we have
$$H^0_{\L}(\O)=I(\O)=C^{\infty}({\R}^{\ast}_{+}).$$
Now, note $\b=dr$ and let $X$ be the Euler vector field ($X=\sum x_i\p_i$).
Consider the tranversal distribution $\nu{\mathcal F}$ where, for each $x$, 
$\nu_{x}^{\ast}{\mathcal F}$ (resp. $\nu_{x}{\mathcal F}$) is  spanned by
$\b_{x}$ (resp. {$X_{x}$).
Recall also that $H^1_{\L}(\O)$ is reduced to the set:
$$\{ A\in {\mathcal V}_{1,0}(\O): \s''(A)=0 \,\hbox{and}\,
\exists B\in {\mathcal V}_{0,1}(\O)\,\hbox{such that}\,\s(A)+\s''(B)=0\}.$$
It is easy to see that the set
$\{A\in {\mathcal V}_{1,0}(\O):\s''(A)=0\}$ coincides with the set $\{ u X: u\in I(\O)\}$.
Now, since $\s(u X)=u \s(X)=u \L$ for all $u$ in $I(\O)$ and that the class 
$[u \L]$ in $H_{\L,tan}^2(\O)$ does not vanish unless $u=0$, we obtain $$H^1_{\L}(\O)=\{0\}.$$\par\noindent
Using Proposition 6 (i), we get
$$H^2_{\L}(\O)\simeq \frac{H^2_{\L,tan}(\O)}
{\s({\mathcal V}_{1,0}(\O))}\oplus H^1(\O,\P^1({\mathcal
F})).$$ We saw that ${\mathcal V}_{1,0}(\O)\supset\{u X: u\in I(\O)\}$ and that $\s(u X)=u \L$ for all
$u$ in $I(\O)$. Moreover, by Proposition 6 (ii), 
$H^1(\O,\P^1({\mathcal F}))$ and $H^1_{\L,tan}(\O)$ are isomorphic. It follows that
$$H^2_{\L}(\O)=\{0\}.$$\par\noindent
Finally, again by Proposition 6, we have
$$H^3_{\L}(\O)\simeq H^2(\O,\P^1({\mathcal F}))\quad\hbox{and}\quad H^2(\O,\P^1({\mathcal F}))\simeq H^2_{\L,tan}(\O).$$
More explicitly, 
\begin{align*}
H^3_{\L}(\O)&=\{u [\L\w X]: u\in I(\O)\}\cr
&\simeq \{u [w_{\L}\w dr]:u\in I(\O)\}\cr
&\simeq I(\O).
\end{align*}
Lastly,
$$H^k_{\L}(\O)=\{0\}\quad\forall k>3.$$\noindent 

\begin{Rem}
The Lie algebra $\frak{su}(2)$ provides an example of a regular Poisson manifold, namely
$\O=\frak{su}(2)^{\ast}\backslash\{0\}$, which is exact ($\L=\s(X)=\s(\sum x_i\p_i))$ without being tangentially
exact ($\L$ cannot be written in the form $\L=\s(T)$ for any tangential vector field $T$).\par\noindent
In fact, since for $\O=\frak{su}(2)^{\ast}\backslash\{0\}$, $H^2_{\L}(\O)=\{0\}$ and $H^2_{\L,tan}(\O)\not=\{0\}$, 
the case of $\frak{su}(2)$ illustrates the fact that the TP-cohomology spaces are generally not imbedded 
in the corresponding P-cohomology spaces, 
except
in degree 1.\par\noindent
Note also that there is in [Xu] a different and beautiful method to calculate the P-cohomology
of $\frak{su}(2)^{\ast}\backslash\{0\}$ by means of symplectic groupoids.
It consists of converting the P-cohomology to the de Rham cohomology of certain manifolds.
\end{Rem}

\medbreak
We propose now to discuss the remaining 3-dimensional Lie algebras. As we shall see in the sequel, all of 
them are pathological cases.\medbreak

\

\subsection{An interesting pathological example}

\

Consider the Lie algebra $\frak{h}$ given by the following brackets:
$$[X_1,X_2]=X_2,\quad [X_1,X_3]=-X_3$$
and denote by $H$ the connected and simply connected Lie group with Lie algebra $\frak{h}$.
For this Lie algebra,
the 2-dimensional coadjoint orbits are the connected components of the hyperbolic cylinders 
$x_2x_3=const.$ and
the half planes
$x_2=0$, $sign(x_3)$ fixed and  $x_3=0$, $sign(x_2)$ fixed. 
Each of them is cohomologically trivial and 
the regular Poisson manifold $(\O,\L)$ associated to $\frak{h}$ is the set
$$\O=\{(x_1,x_2,x_3)\in\frak{h}^{\ast}: x_2^2+x_3^2\not=0\},$$
endowed with its regular Poisson structure $\L$.
However, Theorem 2  cannot be applied because, as it was the case for $\frak{g}_{4,1}$,
the space of leaves $\O/H$ is not Hausdorff. \par\noindent
What is the TP-cohomology of $\O$?\par\noindent 
As always, $H^0_{\L,tan}(\O)=I(\O)$.
To describe $H^1_{\L,tan}(\O)$, we proceed as we did for $\frak{g}_{4,1}$. 
We first observe that every tangential vector field $X$ is of the form
$X=a\p_1 + b  H_{x_1}$ with $a,b$  in 
$C^{\infty}(\O)$ and satisfies the equality $\s (X)=0$ if and only if $H_{x_1}(b)+\p_1(a)=0$.
Then, we fix a tangential vector field $X=a\p_1+b H_{x_1}$ such that $\s (X)=0$.
 $X$ can be written in the form $X=\s f$
$(f\in C^{\infty}(\O))$ if and only if there exists $g(x)=g(x_2,x_3)$ in $C^{\infty}(\O)$ such that
$H_{x_1}(g)(x)=a(0,x_2,x_3)$.\par\noindent
Consider now the change of variables: $u=x_2x_3$, $v=\frac{x_2^2-x_3^2}{2}$.
It gives rise to a diffeomorphism $\varphi$ from
$\{(u,v): u>0\}$ to $\{(x_2,x_3): x_2>0,x_3>0\}$ defined by
\begin{align*}
\v(u,v)&=(\sqrt{v+\sqrt{u^2+v^2}},\frac{u}{\sqrt{v+\sqrt{u^2+v^2}}})\quad \hbox{if}\quad v>0\cr
\v(u,v)&=(\frac{u}{\sqrt{-v+\sqrt{u^2+v^2}}},\sqrt{-v+\sqrt{u^2+v^2}})\quad\hbox{if}\quad v<0\cr
\v(u,0)&=(\sqrt{u},\sqrt{u}).
\end{align*}
 $\v$ can naturally be extended to
$\{(u,v):u>0\}\cup \{(0,v):v\not=0\}$.
We still denote by $\v$ this natural extension.\par\noindent
The following result is, for $\frak{h}$, the analog of Lemma 3 given in Section 3 for $\frak{g}_{4,1}$.
\begin{Lem}
Let $a(x)=a(x_2,x_3)$ be in $C^{\infty}(\O)$. Suppose there exists some
function $g(x)=g(x_2,x_3)$ in
$C^{\infty}(\O)$ such that $H_{x_1}(g)=a$ and note $A=a\circ\v$, where
$\v$ is as above.
Then, the limit
$$\lim_{u\rightarrow 0; u>0}\int_{-1}^{1} \frac{A(u,t)}{2\sqrt{u^2+t^2}} dt,$$ 
exists.
\end{Lem}\par
{\sl Proof.}
By assumption, $H_{x_1}(g)= (x_2\p_2-x_3\p_3)(g)=a$. Changing variables, we get:
$$2\sqrt{u^2+v^2}\p_v(G)=A$$ where $G=g\circ \v$.
Thus, on the set $\{(u,v):u>0\}$, $G$ is necessarily of the form
$$G(u,v)=\int_{-1}^v \frac{A(u,t)}{2\sqrt{u^2+t^2}} dt + \psi(u).$$
Moreover, since $\v(0,-1)=(0,\sqrt{2})$ and $g$ is continuous at $(0,\sqrt{2})$,
$\dis\lim_{u\rightarrow 0; u>0}G(u,-1)$ exists. Thus, 
$\dis\lim_{u\rightarrow 0; u>0} \psi(u)$ exists too. In the same way, 
$\dis\lim_{u\rightarrow 0; u>0} G(u,1)$ exists and  
Lemma 3 is proved.\medbreak
\noindent Now, using the vector fields
$\tilde{T_{\a}}=\tilde{t_{\a}}\p_1$ and
$\tilde{K_{\a}}=\tilde {k_{\a}}\p_1$ where
$$\tilde{t_{\a}}= \frac{1}{(x_2^2+x_3^2)^{\a}}\quad\hbox{and}\quad
 \tilde{k_{\a}}=\frac{\exp(\frac{\a}{(x_2^2+x_3^2)^2})}{(x_2^2+x_3^2)^3},$$
one can see that the space $H^1_{\L,tan}(\O)$ is infinite dimensional.\par\noindent
As for the other TP-cohomology spaces, $$H^k_{\L,tan}(\O)=\{0\} \quad\forall k>1.$$
\smallbreak\noindent
To compute the P-cohomology of $\O$, let us say that the symplectic foliation ${\mathcal F}$ of $\O$
is orientable and is defined by the nonsingular 1-form $\b=x_2dx_3+x_3dx_2$ ($d\b=0$).
By Proposition 6 and if 
we note $X=\dis\frac{x_2 \p_3+ x_3\p_2}{x_2^2+x_3^2}$, we have
\begin{align*}
H^0_{\L}(\O)&=I(\O)\cr
H^1_{\L}(\O)&= H_{\L,tan}^1(\O) \oplus \{ u X: u\in I(\O)\}\cr
H^2_{\L}(\O)&\simeq  H^1(\O,\P^1({\mathcal F}))\simeq H^1_{\L,tan}(\O)\cr
H^3_{\L}(\O)&\simeq H^2(\O,\P^1({\mathcal F}))\simeq H^2_{\L,tan}(\O)=\{0\}\cr
H^k_{\L}(\O)&= \{0\}\quad \forall k> 3.
\end{align*}
Thus, the large space $H_{\L,tan}^1(\O)$ happens in the P-cohomology. 
In other words, the spaces $H_{\L}^1(\O)$ and $H_{\L}^2(\O)$ are infinite dimensional.  \smallbreak\noindent
\begin{Rem}
Consider the family of 3-dimensional Lie algebras
given by:
$$[X_1,X_2]=X_2,\quad [X_1,X_3]=\frac{1}{\tau} X_3\quad\tau <-1.$$
One can prove that the TP-cohomology of $\O$ in these examples is essentially the same as in the
case of $\frak{h}$ because of a similar leaf structure.
The situation seems more complicated for the P-cohomology since Proposition 6 cannot be used.
\end{Rem}

\subsection{The case of $\frak{sl}(2)$}

\

We want now to show how the  TP-cohomology computation for
$\frak{sl}(2)$ leads to the same conclusions as in the case of $\frak{g}_{4,1}$ or $\frak{h}$.\par
The Lie algebra $\frak{sl}(2)$ is given by the brackets:
$$[X_1,X_2]=- X_3,\quad [X_3,X_1]=X_2,\quad [X_3,X_2]=- X_1.$$
Let us  identify $\frak{sl}(2)$  with its dual $\frak{sl}(2)^{\ast}$ by means of the Killing form.
Then the orbit decomposition is as follows. There are three orbits in the light cone:  
\begin{align*}
W_0&= \{0\};\cr
W_{+}&=\{(\mu_1,\mu_2,\mu_3): \mu_1^2+\mu_2^2-\mu_3^2=0, \mu_3>0\};\cr
W_{-}&=\{(\mu_1,\mu_2,\mu_3): \mu_1^2+\mu_2^2-\mu_3^2=0, \mu_3<0\}.
\end{align*}
Moreover, the hyperbolic orbits are single sheeted hyperboloids outside the cone:
$$W_{k}=\{(\mu_1,\mu_2,\mu_3): \mu_1^2+\mu_2^2-\mu_3^2=k^2\}\quad (k>0)$$
and double sheeted hyperboloids inside the cone
(each of the sheets being a different orbit):
$$W_{h}=\{(\mu_1,\mu_2,\mu_3): \mu_1^2+\mu_2^2-\mu_3^2=-h^2, \mu_3 h>0\}\quad (h\not=0).$$
The regular Poisson manifold $(\O,\L)$ associated to $\frak{sl}(2)$ is thus the set
$$\O=\{(x_1,x_2,x_3)\in \frak{sl}(2)^{\ast}:x_1^2+x_2^2+x_3^2\not=0\},$$
with its regular Poisson structure $\L$ and the space of leaves is clearly not Hausdorff.  
\par\noindent
As always, $H^0_{\L,tan}(\O)=I(\O)$.\par\noindent
How about $H^1_{\L,tan}(\O)$?
Assume first that $X$  is a tangential vector field of the form
$$X=a H_{x_3} $$
where $a$ is in $C^{\infty}(\O)$.
One can easily check that $X=\s f$ for a function $f$ in $C^{\infty}(\O)$ if and only if  $f$
satisfies the following properties:
\begin{align*}
H_{x_1}(f) &= x_2 a\cr
H_{x_2}(f) &= -x_1 a\cr
H_{x_3}(f) &= 0.
\end{align*}
Consider the open set $U=\{(x_1,x_2,x_3)\in \O: x_1^2+x_2^2-x_3^2>0\}$ which can be parameterized by:
\begin{align*}
x_1 &=-p \sin (q)+z \cos (q)\cr
x_2 &=-p\cos(q)-z \sin(q)\cr
x_3 &=p.
\end{align*}
Using this parameterization, we see that, if $X=\s f$, then 
$\p_{p}(f)=- a$  and $\p_{q}(f)=0$.
In other words, if $X=\s f$, then $f_{|U}$ must depend only on $(p,z)$ and be of the form
$$f=f(p,z)=-\int_{-1}^p a(s,q,z) ds + \psi(z).$$
Consider the vector fields 
$\hat{T_{\a}}$
ans $\hat{K_{\a}}$, $\a>\frac{1}{2}$, defined by
$$\hat{T_{\a}}= \hat{t_{\a}} H_{x_3}\quad\hbox{and}\quad 
\hat{K_{\a}}= \hat{k_{\a}} H_{x_3}$$ where
$$\hat{t_{\a}}= \frac{1}{(x_1^2+x_2^2+x_3^2)^{\a}}\quad\hbox{and}\quad 
\hat{k_{\a}}=\frac{\exp(\frac{\a}{x_1^2+x_2^2+x_3^2})}{(x_1^2+x_2^2+x_3^2)^2}.$$
Then, with the same arguments as in the examples of $\frak{g}_{4,1}$ and $\frak{h}$,
one can prove that the space $H_{\L,tan}^1(\O)$ is infinite dimensional.
\medbreak\noindent
The same is true for $H^2_{\L,tan}(\O)$.
Indeed, take  a tangential 2-tensor field $A$. Necessarily, there is a function $\v$ in $C^{\infty}(\O)$ such that
$$A=\v \L=\v (x_3 \p_2\w\p_1+ x_2\p_3\w\p_1+x_1\p_2\w\p_3)$$
and $A$ satisfies  $\s(A)=0$.
Now, suppose there exists  a tangential vector field $B$ 
such that $A=\s (B)$. If we note $B=a H_{x_1}+b H_{x_2}+ c H_{x_3}$ where $a,b,c$ are in $C^{\infty}(\O)$, 
then we obtain
$$(\#)\quad H_{x_1}(a)+H_{x_2}(b)+H_{x_3}(c)=\v.$$
Translating $(\#)$ on the open set  $U=\{(x_1,x_2,x_3)\in \O:x_1^2+x_2^2-x_3^2>0\}$, we get
$$(-x_2\p_p-\sin (q) \p_q)(a)+ (x_1\p_p-\cos (q)\p_q)(b)+\p_q(c)=\v.$$
Let us integrate the above expression  between 0 and $2\pi$ as follows
\begin{align*}
&p\p_p \Bigl(\int_0^{2\pi} (\cos (q) a-\sin (q) b) dq\Bigr) + z \p_p \Bigl(\int_0^{2\pi} (\sin (q) a  
+\cos(q) b) dq\Bigr)\cr
&\quad + \quad \int_0^{2\pi} (-\sin (q) \p_q(a)-\cos (q)\p_q(b)) dq=\cr
&= p\p_p \Bigl(\int_0^{2\pi} (\cos (q) a -\sin (q) b) dq\Bigr) + z \p_p \Bigl(\int_0^{2\pi} (\sin (q) a 
 +\cos (q) b) dq\Bigr)\cr
&\quad + \quad \int_0^{2\pi} (\cos (q) a- \sin (q) b) dq=\cr
&= \int_0^{2\pi} \v(p,q,z) dq.
\end{align*}
Assume that $\v$ depends only on the variables $p,z$.
Moreover, define the functions $H$ and $g$ by 
$$H(p,z)=\dis\int_0^{2\pi} (\cos (q) a-\sin (q) b) dq,\quad 
g(p,z)=\dis\int_0^{2\pi} (\sin (q) a+\cos (q) b) dq.$$
Then we have
$$ p\p_p H+H= 2\pi \v -z \p_p(g).$$
Thus, on the open set $U\cap \{p>0\}$, $H$ is necessarily  of the form
$$H= \frac{1}{p} \Bigl(\int_{1}^p 2\pi \v(s,z) ds - z g(p,z) +z g(1,z)+\psi_1(z)\Bigr).$$
Similarly, on the  open set $U\cap \{p<0\}$, $H$ is necessarily of the form
$$H= \frac{1}{p} \Bigl(\int_{-1}^p 2\pi \v(s,z) ds - z g(p,z) +z g(-1,z)+\psi_2(z)\Bigr).$$
The function $H(p,z)$ can be extended at the points $(p,0)$ with $p\not=0$,
therefore $\psi_1(z)$ et $\psi_2(z)$
have a limit when $z(>0)$ tends to zero.
Moreover, $H(p,z)$ can also be extended at the points $(0,z)$ with $z\not=0$.
Thus, for all $z\not=0$, we get
$$\int_1^0 ( 2\pi \v(s,z)) ds - z g(0,z) +z g(1,z)+\psi_1(z)=0$$
and
$$\int_{-1}^0 ( 2\pi \v(s,z)) ds - z g(0,z) +z g(-1,z)+\psi_2(z)=0.$$
It follows that
$$\int_{-1}^1 (2\pi \v(s,z)) ds - zg(1,z)+z g(-1,z)-\psi_1(z)+\psi_2(z)=0.$$
In other words, if $A=\v \L$ where $\v_{|U}=\v(p,z)$ and if there exists a tangential vector field $B$ such that
$A=\s (B)$, then the limit
$$\lim_{z\rightarrow 0; z>0} \int_{-1}^1 2\pi \v(s,z) ds$$ should exist.
Now to see that the space $H_{\L,tan}^2(\O)$ is infinite dimensional, 
one can take the vector fields 
$\hat{t_{\a}} \L$ where $\hat{t_{\a}}=\dis\frac{1}{(x_1^2+x_2^2+x_3^2)^{\a}}$
and  $\hat{k_{\a}}\L$ where
$\hat{k_{\a}}=\dis\frac{\exp(\frac{\a}{ x_1^2+x_2^2+x_3^2})}{(x_1^2+x_2^2+x_3^2)^2}.$
\par\noindent
Lastly, $$H_{\L,tan}^k(\O)=0\quad\forall k>2.$$
Let us also say that the symplectic foliation ${\mathcal F}$ of $\O$ is orientable and is defined
by the nonsingular 1-form $\b=x_1 dx_1+x_2dx_2-x_3dx_3$, which satisfies 
$d\b=0$.
Thus, using Proposition 6, one can see that
all the P-cohomology spaces $H_{\L}^k(\O)$ ($1\leq k\leq 3)$ are infinite dimensional.
\par\smallskip
We finish with an application of our cohomology calculations for the classification of tangential star products:

\smallskip\noindent                                                                                      
\begin{Rem}
As we just proved in this section, the second TP-cohomology space of $\O$ is zero                                                                                       
for any 3-dimensional Lie algebra except for $\frak{su}(2)$ and $\frak{sl}(2)$. 
This implies the uniqueness (up to equivalence) of the
tangential star products on $\O$ for any non simple 3-dimensional Lie algebra.
\end{Rem}

\

\section{Concluding remarks}\label{sec5}

As shown in Section 2, for a regular Poisson manifold, the TP-cohomology coincides
with the leafwise de Rham cohomology of the symplectic foliation, thus
unlike the P-cohomology, does not depend on the symplectic structure
along the leaves.
That relates the task of computing the TP-cohomology of regular Poisson manifolds to non-trivial
questions of foliation theory.
Roughly speaking, the TP-cohomology not only contains the de Rham cohomology
of the leaves, but also translates the foliation complexity which
includes essentially the relative position of the leaves.
As seen from the computations of Sections 3 and 4, when the space of leaves is not Hausdorff, this TP-cohomology is
very large and hardly describable.
To finish, let us say that the TP-cohomology spaces are involved in the P-cohomology,
for instance we have the inclusion $H^1_{\L,tan}(M)\subset H^1_{\L}(M)$ for any
regular Poisson manifold $M$. Thus,
our calculations and comments can be of some help to understand why the P-cohomology itself is,
as often said, so difficult to compute.

\

\vfill\eject
\centerline{{\bf Acknowledgements.}}\vskip0,25cm\noindent
I am infinitely grateful to Didier ARNAL for his continual and helpful guidance throughout 
this work. I thank also the referee of this paper for illuminating comments.
The research has been supported by an ARC of the Communaut\'e Fran\c{c}aise de Belgique.

\

\centerline{{\bf
References.}}\vskip0,25cm\noindent [APP] J.A. de Azcarraga, A.M. Perelomov and J.C. Perez Bueno, 
{\sl The Schouten Nijenhuis bracket, cohomology and generalized
 Poisson structures}, J. Phys. A.:
Math. Gen.,
 {\bf 29} (1996), 7993-8009.\par\noindent
[ACG] D. Arnal, M. Cahen and S. Gutt, {\sl Deformations on coadjoint orbits}, J. Geom. Phys., 
{\bf vol.3}, {\bf 3} (1986), 327-351.\par\noindent
[BLM] J. Boidol, J. Ludwig and D. M\"uller, {\sl On infinitely small orbits}, Studia Mathematica, T. LXXXVIII (1988).
\par\noindent
[Bon] P. Bonnet, {\sl Param\'etrisation du dual d'une alg\`ebre de Lie nilpotente},
Ann. Inst. Fourier Grenoble, {\bf 38}, {\bf 3} (1988), 169-197.\par\noindent
[Br]  R. Bryant, {\sl An introduction to Lie groups and symplectic geometry}, Freed, Daniel S. (ed.) et al., Geometry and 
quantum field theory. Lecture notes from the graduate summer school program,
June 22-July 20, 1991, Park City, UT, USA. Providence, RI:American Mathematical Society. 
IAS/Park City Math.Ser.1 (1995), 5-181.
\par\noindent
[Bry] J.L. Brylinski, {\sl Loop Spaces, Characteristic Classes and Geometric Quantization}, Progress in Mathematics,
{\bf  vol.107}, Birkh\"auser, Boston, 1993.\par\noindent
[CGR] M. Cahen, S. Gutt and J. Rawnsley, {\sl On tangential star products for the coadjoint Poisson structure},
Comm. Math. Phys., {\bf 180} (1996), 99-108.\par\noindent
[CW] A. Cannas Da Silva and A. Weinstein, {\sl Geometric Models for Noncommutative Algebras}, Berkeley
mathematics lecture notes, {\bf vol.10}, American Mathematical Society (1999).\par\noindent
[CG] L. Corwin and F.P. Greenleaf, {\sl Representations of nilpotent Lie groups and their applications.}
Part 1: Basic theory and examples. Cambridge studies in advanced mathematics, 1990.\par\noindent
[DH] P. Dazord et G. Hector, {\sl Int\'egration Symplectique des vari\'et\'es de Poisson totalement
 asph\'eriques.} In: 
Symplectic Geometry, Groupoids and Integrable Systems, S\'eminaire Sud Rhodanien 
de G\'eom\'etrie \`a Berkeley
(1989) (P. Dazord and A. Weinstein, eds.). MSRI Publ.20, Springer-Verlag, Berlin-Heidelberg-New York (1991), 
37-72.
\par\noindent
[Die] J. Dieudonn\'e, {\sl El\'ements d'analyse}, Cahiers scientifiques, fasc.
 XXXIII, tome 3, Gauthier-Villars, 1974.
\par\noindent
[GK] M. Goze and Y. Khakimdjanov, {\sl Nilpotent Lie algebras}, Mathematics and its Applications, {\bf 361},
Kluwer Academic Publishers Group, 1996.\par\noindent
[Gr\'e] G. Gr\'elaud, {\sl Sur les repr\'esentations des groupes de Lie r\'esolubles}, Th\`ese,
Universit\'e de Poitiers (1984).\par\noindent
[Hue] J. Huebschmann, {\sl Poisson cohomology and quantization}, J. f\"ur Reine und Angew. Math., {\bf 408} (1990), 57-113.
\par\noindent
[Kon] M. Kontsevich, {\sl Deformation quantization of Poisson manifolds}, Preprint q-alg/9709040 (1997).\par\noindent
[Kos] J.L. Koszul, {\sl Crochet de Schouten-Nijenhuis et cohomologie},  Elie Cartan et les 
Math\'ematiques d'aujourd'hui, 
Ast\'erisque hors s\'erie, Soci\'et\'e Math. de France (1985), 257-271.\par\noindent
[Li1] A. Lichnerowicz, {\sl Les vari\'et\'es de Poisson et leurs alg\`ebres de Lie associ\'ees}, 
 J. Diff. Geom., {\bf 12} (1977), 253-300.\par\noindent
[Li2] A. Lichnerowicz, {\sl Vari\'et\'es de Poisson et feuilletages}, Ann. Fac. Toulouse (1982), 195-262.
\par\noindent
[Mas]
M. Masmoudi, {\sl Tangential formal deformations of the Poisson bracket and tangential star products on a regular 
Poisson manifold}, J. Geom. Phys., {\bf 9} (1992), 155-171.\par\noindent
[Pe1] N.V. Pedersen, {\sl 
Geometric quantization and the universal enveloping algebra of a nilpotent Lie group},
Transactions of the American Mathematical Society, {\bf 315}, {\bf 2} (1989), 511-563.\par\noindent
[Pe2] N.V. Pedersen, {\sl Geometric quantization and nilpotent Lie groups}, A collection of 
examples, University of Copenhagen Denmark (1988), 1-180.\par\noindent
[Puk] L. Pukanszky, {\sl Le\c{c}on sur les repr\'esentations des groupes}, Dunod, Paris,
1967.\par\noindent
 [SG] M. Saint-Germain, {\sl Poisson algebras and transverse structures},
 J. Geom. Phys., {\bf 31} (1999), 153-194.\par\noindent
[Va1] I. Vaisman, {\sl Remarks on the Lichnerowicz-Poisson cohomology},  Ann. Inst. Fourier Grenoble, {\bf 40} (1990),
 951-963. \par\noindent
[Va2] I. Vaisman, {\sl Lectures on the Geometry of Poisson Manifolds}, Progress in Mathematics,
{\bf vol.118}, Birkh\"auser, Basel, 1994.\par\noindent
[Ver] M. Vergne, {\sl La structure de Poisson sur l'alg\`ebre sym\'etrique d'une alg\`ebre de Lie nilpotente}, 
Bull. Soc. Math. Fr., {\bf 100} (1972), 301-335.\par\noindent
[Wei] A. Weinstein, {\sl The local structure of Poisson manifolds}, J. Diff. Geom., {\bf 18} (1983),
523-557.\par\noindent [Xu] P. Xu, {\sl Poisson cohomology of regular Poisson manifolds}, Ann. Inst.
Fourier Grenoble, {\bf 42} (1992),
 967-988.

\end{document}